\documentclass[12pt]{article}
\usepackage{amsmath}
\usepackage{amsfonts}
\usepackage{amssymb}
\usepackage{citesort}
\usepackage{graphicx}
\usepackage{subfigure}%

\begin{document}

\newtheorem{theorem}{Theorem}
\newtheorem{acknowledgement}[theorem]{Acknowledgement}
\newtheorem{algorithm}[theorem]{Algorithm}
\newtheorem{problem}[theorem]{Problem}

\title
{Fast reconstruction algorithms for the thermoacoustic tomography
in certain domains with cylindrical or spherical symmetries}

\author{Leonid Kunyansky}
%\address{Department of Mathematics \\
%University of Arizona\\
%Tucson, AZ 85721, USA}
%\email{leonk@math.arizona.edu}

%\subjclass{Primary: 44A12, 65M32; Secondary: 92C55}
%\keywords{Radon transform, spherical means, thermoacoustic tomography,
%fast algorithms, integrating detectors}

\maketitle

\begin{abstract}
We propose three fast algorithms for solving the inverse problem of the
thermoacoustic tomography corresponding to certain acquisition geometries.
Two of these methods are designed to process the measurements done
with point-like detectors placed on a circle (in 2D) or a sphere (in 3D)
surrounding the object of interest. The third inversion algorithm
works with the data measured by the integrating line detectors arranged
in a cylindrical assembly rotating around the object.
The number of operations required by these techniques is equal
to
$\mathcal{O}(n^{3} \log n)$ and
$\mathcal{O}(n^{3} \log^2 n)$ for the 3D techniques
(assuming the reconstruction grid with $n^3$ nodes) and
to $\mathcal{O}(n^{2} \log n)$ for the 2D problem with $n \times n$
discretizetion grid.
Numerical simulations show that our methods are at least two orders of magnitude
faster than the existing algorithms, without any sacrifice in accuracy or
stability. The results of reconstructions from real measurements
done by the integrating line detectors are also presented,
to demonstrate the practicality of our methods.

\end{abstract}

\section*{Introduction}

Thermoacoustic tomography (TAT) and the closely related photoacoustic modality
(PAT) are both based on the measurements of acoustic waves generated within the
object of interest by a thermoelastic expansion \cite{Oraev94,Kruger}. To
initiate the wave, the object is illuminated by a short electromagnetic (EM)
pulse whose energy is partially absorbed by the tissues.
The subsequent increase in the local temperature makes the tissues expand,
which in turn generates the outgoing acoustic wave registered by the detectors
surrounding the object. By solving the corresponding inverse problem one can
reconstruct the distribution of the initial pressure inside the body.
Spikes in the initial pressure are indicative of cancerous tumors that
absorb much more EM energy than healthy tissues; thus, cancer detection
(particularly, in breast imaging) is one of the most promising
applications of this modality.
General review of the inverse problem of TAT/PAT, and investigation of
the related mathematical questions can be found in reviews
\cite{CRC,Wang-book,Springer} and references therein.

The recent advances in the measuring technology of TAT/PAT make it possible to collect
enough data to reconstruct high resolution 3D images of the object of
interest. Since such a reconstruction amounts to computing many millions
of unknowns, the use of existing reconstruction algorithms may lead
to inordinate computation time. For example, methods based of
filtration/backprojection formulas would require several days of computations
per one high-resolution 3D image. We, thus, propose three fast reconstruction
algorithms for certain acquisition geometries with spherical or circular
symmetry, including, in particular, a fast method for the measuring scheme
with integrating line detectors arranged in a rotating cylindrical assembly.
The present algorithms produce high quality images two orders of magnitude
faster than the existing methods.

Let us briefly review existing reconstruction techniques while paying
close attention to the asymptotic estimates
of the number of floating point operations (flops) needed to complete one
reconstruction. (While such estimates hide the unknown constant factor,
they usually provide a good qualitative measure of the efficiency of the method.)
Ideally, one would like to reconstruct an image on a grid with $N$
unknowns  in $\mathcal{O}(N)$ flops. However, methods that
require $\mathcal{O}(N\log^{\alpha}N)$ flops (where $\alpha$ is some constant)
are still generally considered "fast".

The simplest inverse problem in TAT/PAT arises when small (point-like)
detectors are placed on an infinite plane. In this case the explicit solution
can be obtained using the Fourier transform techniques (see
\cite{Andersson,Fawcett,Nat-Wub} and references therein). The authors of
\cite{Haltm-fft}, by combining proper discretization of such a solution with
application of the non-uniform FFT (NUFFT) developed a fast algorithm that
reconstructs images on $n\times n\times n$ grid (in 3D) in $\mathcal{O}%
(n^{3}\log n)$ flops. However, in any practical application such a measuring
plane needs to be truncated, which leads to a loss of information about the
wave fronts nearly parallel to the plane. Thus, it is
preferable to use closed (and bounded) measuring surfaces.

One of the simplest closed surfaces is a sphere, and first explicit solutions
of the TAT/PAT\ problem in a closed domain were obtained for circular (in 2D)
and spherical (in 3D) acquisition geometries in \cite{Norton1} and
\cite{Norton2} by means of separation of variables.
Later, some modifications
of the 2D formulas of \cite{Norton1} were proposed in \cite{AmbKuch} and
\cite{Haltm-circ} in order to avoid the divisions by zero in the original
formulas of \cite{Norton1}.
% the details of the corresponding reconstruction
%algorithm and the numerical results were given in \cite{Haltm-circ}.
Such
techniques have complexity $\mathcal{O}(n^{3})$ flops for a 2D grid with
$n^{2}$ unknowns.
%The 3D version of the series solution \cite{Norton2} was not
%formulated as an algorithm, but its straightforward implementation would
%result in a rather slow $\mathcal{O}(n^{6})$ algorithm. However, a
%modification of this technique proposed in \cite{Ramm}, if properly
%discretized, can lead to an $\mathcal{O}(n^{4})$ algorithm.
A straightforward implementation of the 3D version of the series solution \cite{Norton2}
results in a rather slow $\mathcal{O}(n^{6})$ algorithm. However, a
modification of this technique proposed in \cite{Ramm}, if properly
discretized, yields a faster $\mathcal{O}(n^{4})$ algorithm.

There also exist several explicit inversion formulas of filtration/backprojection type
\cite{FPR,Finch-even,Kunyansky,MXW2,nguyen,burg-FBP}
for the  spherical or circular acquisition geometries.
In spite of the
theoretical importance of these formulas, the corresponding algorithms
are not fast, requiring
$\mathcal{O}(n^{3})$ flops for a 2D reconstruction on
$n\times n$ grid, and $\mathcal{O}(n^{5})$ flops for a 3D problem.

The time reversal by means of a finite difference solution of the wave equation
back in time \cite{HKN,burg-exac-appro} is faster (at least in 3D).
This technique can be adapted to almost any closed
acquisition surface, and it allows one to take into account (known) variations
in the speed of sound within the region of interest (most other
techniques are applicable only if the speed of sound is a known constant).
These methods have complexity $\mathcal{O}(n^{4})$ in 3D and $\mathcal{O}%
(n^{3})$ in 2D.

Mathematically equivalent to the time reversal are the methods based on
expanding the solution of the wave equation into the series of eigenfunctions
of the Dirichlet Laplacian in the domain surrounded by the acquisition surface
\cite{Kun-series,AK}.
%The performance of these techniques strongly depends on
%the explicit knowledge of the eigenfunctions and the existence of efficient
%algorithms for the summation of these series.
This technique is computationally efficient only if there exists a fast method
for the summation of the eigenfunctions. For example, if the detectors
are placed on a surface of a cube, such a summation can be performed using
the 3D Fast Fourier Transform (FFT) algorithm, and one obtains a very
fast $\mathcal{O} (n^{3}\log n)$ reconstruction technique \cite{Kun-series}.
%(this algorithm can
%be modified to solve the 2D problem with the detectors placed on the sides of
%a square, with $\mathcal{O}(n^{2}\log n)$ operation count).
This is the only known fast method for TAT/PAT problems with point-like
detectors placed on a closed surface.

In addition to point-like detectors, there exists another new and promising
class of measuring devices, the so-called integrating line detectors
\cite{burg-area-line,burg-fabry,Palt-Machzend}. A sensing element in a line
detector is a straight optical fiber coupled to a laser interferometer. (In
some implementations the fiber is replaced by a laser beam propagating
through the water in which the object of interest is immersed.) Propagating
acoustic waves slightly elongate the fiber, and this elongation is registered by
the interferometer. The measured value is proportional to the line integral of the
acoustic pressure. Since the fiber can be made very thin, the use of such
detectors can significantly increase the spatial resolution of TAT/PAT.
However, in order to enable this new acquisition technique one needs
new inversion algorithms.

\begin{figure}[t]
\begin{center}
\subfigure[]{\includegraphics[width=2.0in,height=2.0in]{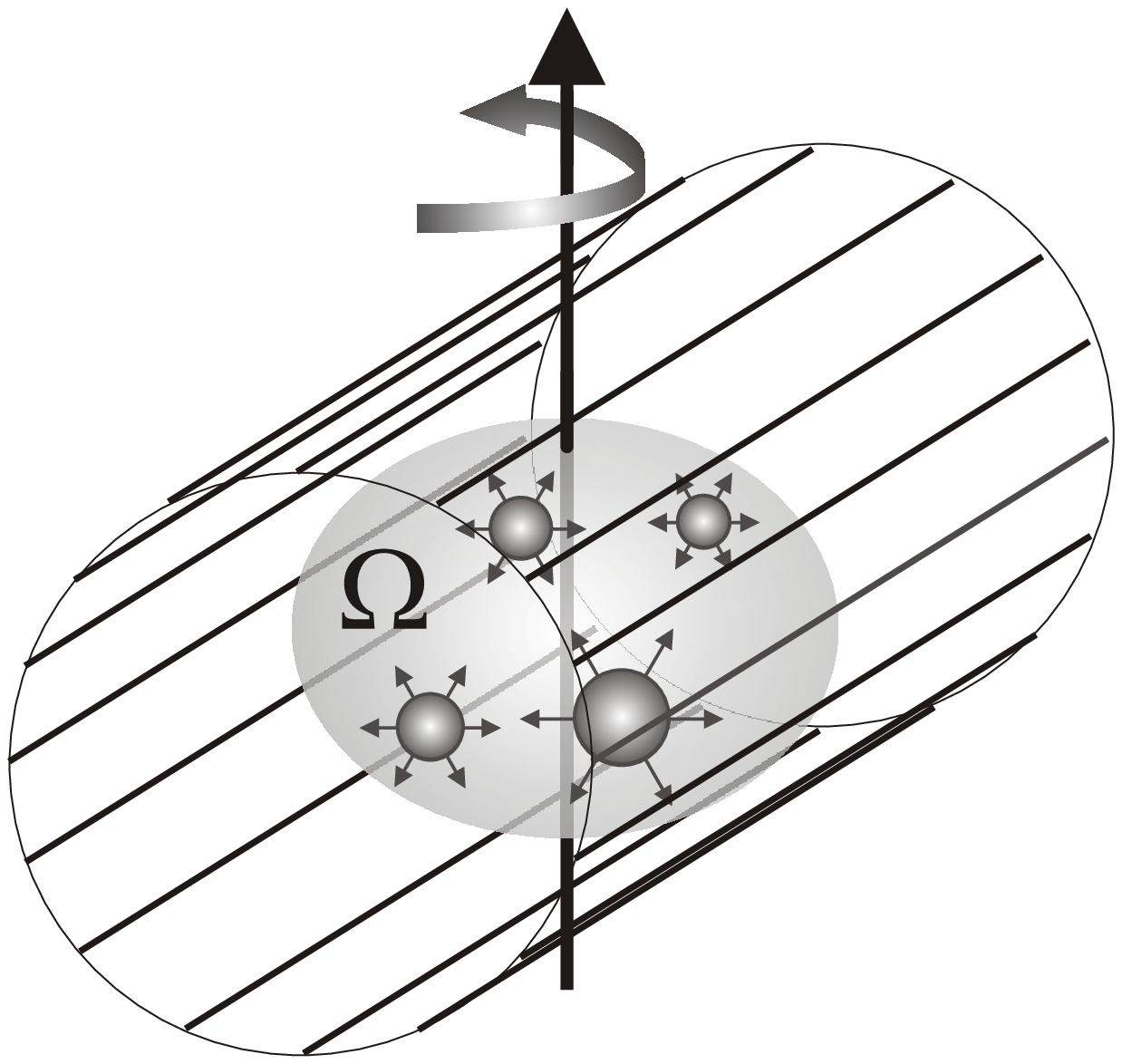} }
\subfigure[]{\includegraphics[width=2.0in,height=2.0in]{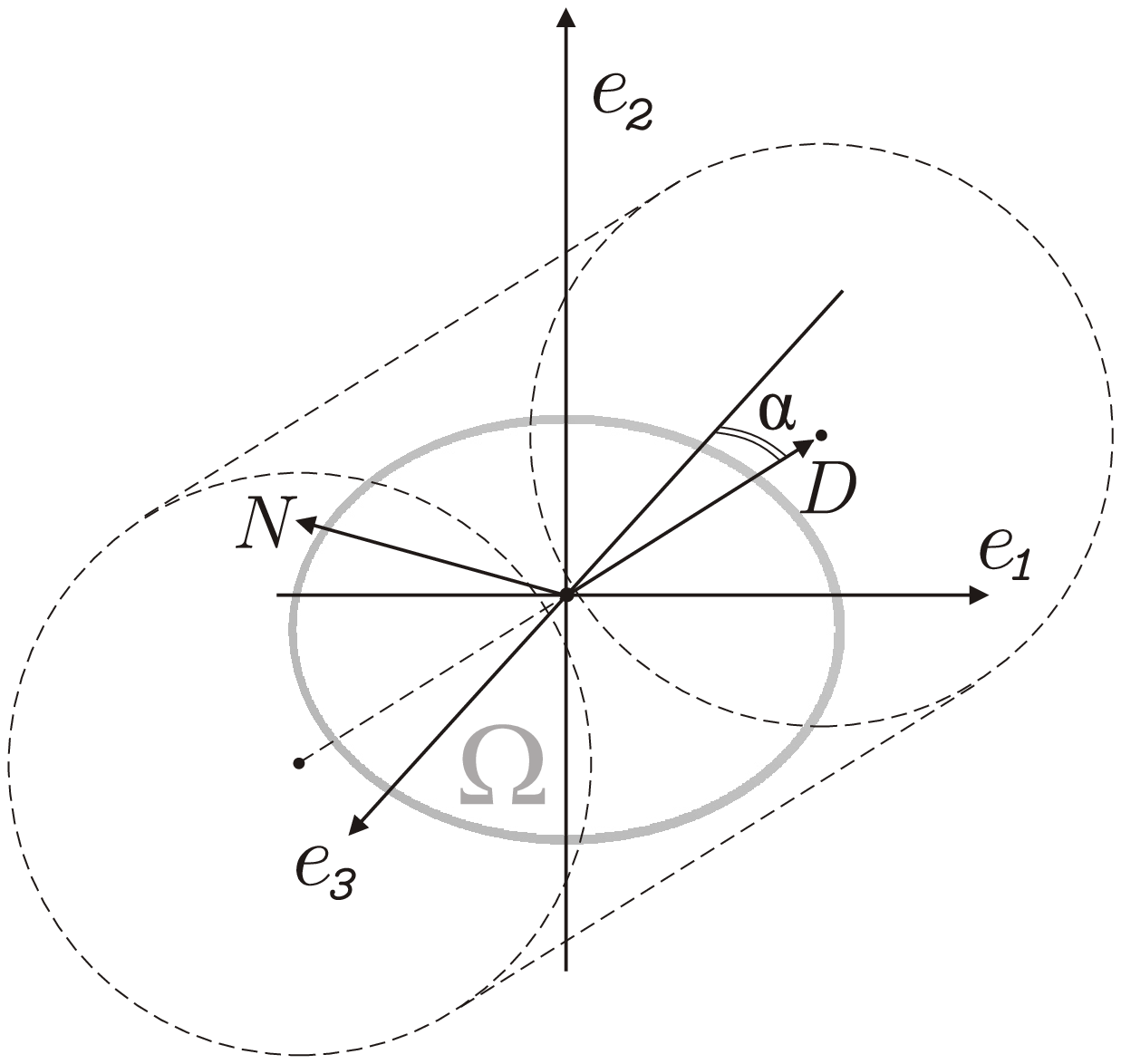}}
\end{center}
\caption{ Measuring scheme with linear integrating detectors (a) general view
(b) geometry}%
\label{F:geom}%
\end{figure}

The most studied measurement scheme of this sort
\cite{burg-area-line,burg-fabry,Palt-Machzend} utilizes integrating line
detectors arranged in a rotating cylindrical assembly (see
Figure~\ref{F:geom}(a)), with the object placed inside the cylinder. A two
step inversion procedure was proposed in \cite{Palt-Machzend} for solving the
arising inverse problem (see also
\cite{burg-FBP,burg-exac-appro,Paltaufexp,Paltaufwei}). It is based on the
observation that, for a fixed orientation of the cylinder, the data measured
by the line detectors are related to the line integrals of the sought function
(over the parallel lines) by the two-dimensional wave equation. Thus, the
latter line integrals can be reconstructed by one of the 2D methods for
TAT/PAT\ with conventional point-like detectors placed on a circle. The first
step of the reconstruction procedure consists in solving such a 2D problem for
each orientation of the cylinder. On the second step one performs a series of
inversions of the 2D Radon transform to reconstruct the sought function from
the line integrals whose values were found on the first step. Since the
fastest known algorithms for the circular geometry require at least
$\mathcal{O}(n^{3})$ operations, and the number of problems to be solved is
$\mathcal{O}(n),$ the first step needs at least $O(n^{4})$ flops. The second
step as described in
\cite{Palt-Machzend,burg-FBP,burg-exac-appro,Paltaufexp,Paltaufwei} is also an
$O(n^{4})$ flops procedure; however, since methods for fast ($\mathcal{O}%
(n^{2}\log n))$ inversion of the 2D Radon transform are well-known (see
\cite{Natterer} and references therein), it can be accelerated to
$\mathcal{O}(n^{3}\log n)$ flops. Thus, the total number of operations
($O(n^{4}))$ for the whole inversion technique is determined by the first step.

One of our goals is to develop a fast
algorithm for the data acquisition scheme with integrating line detectors
described in the previous paragraph. To this end, we first
develop a fast $\mathcal{O}(n^{2}\log n)$ reconstruction algorithm for the 2D
problem with point-like detectors located on a circle. Since such a technique
(and its 3D generalization) are of independent interest for problems with
conventional detectors, we present them separately in
Section~\ref{S:conventional}. Further analysis of the 3D problem with the line
detectors shows that the two-step reconstruction procedure mentioned in the previous paragraph
is somewhat suboptimal. Instead, a direct, efficient $\mathcal{O}(n^{3}\log n)$
reconstruction algorithm can be built by modifying the 2D method described in
Section~\ref{S:2Dconventional}; we present this technique in
Section~\ref{S:linedet}.

The results of numerical simulations show that our methods yield very fast
reconstructions without any sacrifice in stability or in the resolution of the
images. In addition to simulations, in Section~\ref{S:real} we illustrate
the work of our algorithms by reconstructing images from the data of real
measurements performed in RECENDT
(Research Center for Non-Desctructive Testing, Linz, Austria).

\section{Formulation of the problems}

We will assume throughout the paper that the speed of sound in the tissue is
constant (in this case one can set it to be equal 1 without loss of
generality), and that the attenuation is negligible. This simplified model is
acceptable in such important applications as breast imaging, and most of the
formulas and algorithms mentioned in the Introduction are based on these
assumptions. (In the situations when the variations in the speed of sound need
to be taken into account, other techniques (e.g., time reversal) should be used.)

Under the above assumptions the acoustic pressure $u(x,t)$ solves the
following initial problem for the wave equation in the whole space%
\begin{equation}
\left\{
\begin{array}
[c]{c}%
u_{tt}=\Delta u,\text{ }x\in\mathbb{R}^{d},\text{ }t\in\lbrack0,\infty),\\
u(x,0)=f(x),\\
u_{t}(x,0)=0.
\end{array}
\right.  \label{IVP}%
\end{equation}
were the initial pressure of the acoustic wave $f(x)$ is the function we seek,
and dimension $d$ equals 2 or 3, depending on the context. We will assume that
$f(x)$ is finitely supported within the region of interest $\Omega,$ and that
the measurements of $u(x,t)$ are done outside $\Omega.$

\subsection{Circular geometry with conventional detectors}

As explained in the Introduction, we start with the problem that arises when
the pressure is measured by the conventional point-like detectors placed on a
circle (or, in 3D, a sphere) $\partial B$ of radius $R,$ and $\Omega$
coincides with the corresponding disk (or ball) $B$. In other words, the data
$P(y,t)$ are defined as follows%
\[
P(y,t)=u(y,t)\left.  {}\right\vert _{y\in\partial B},t\in\lbrack0,\infty).
\]

\begin{problem}
\textrm{\textbf{(2D)}} Reconstruct $f(x)$ from $P\left(  y,t\right)  $,
$x\in\mathbb{R}^{2}$, $y\in\partial B,$ $t\in\lbrack0,\infty).$
\end{problem}
The solution of this problem will serve as an important step in solving problem 2.

We will also present a fast algorithm for solving the 3D version of this problem:
\setcounter{theorem}{0}

\begin{problem}
\textrm{\textbf{(3D)}} Reconstruct $f(x)$ from $P\left(  y,t\right)  $,
$x\in\mathbb{R}^{3}$, $y\in\partial B,$ $t\in\lbrack0,2R].$
\end{problem}

The difference in the observation time used in 2D and 3D versions arises since
in 3D, due to the Huygens principle, $u(x,t)=0$ on $\partial B$ for all
$t>2R.$ We present fast algorithms for solving Problem 1 in 2D and 3D in
Section~\ref{S:conventional}.

\subsection{Acquisition scheme with the line detectors}

Suppose that the acoustic pressure $u(x,t)$ satisfying conditions (\ref{IVP})
is measured by the integrating line detectors lying on the surface of a
cylinder of radius $R$ whose axis is parallel to the vector $D(\alpha
)=(\cos\alpha,0,\sin\alpha)$ and passes through the origin (see
Figure~\ref{F:geom}). If we denote the left normal to $D$ lying in the
horizontal plane (spanned by the vectors $e_{1}=(1,0,0)$ and $e_{3}=(0,0,1)$ )
by $N(\alpha)=(-\sin\alpha,0,\cos\alpha),$ then each detector lies on a line
$l(\alpha,\beta)$ passing through the point $A(\alpha,\beta)=R\cos\beta
e_{2}+R\sin\beta N(\alpha),$ where $e_{2}=(0,1,0)$ is the vertical unit
vector. This detector measures the value proportional to the line integral
$P_{\alpha}(y(\beta),t)$ of $u:$%
\begin{align}
P_{\alpha}(y(\beta),t)  &  =\int\limits_{\mathbb{R}^{1}}u(y_{1}(\beta
)e_{2}+y_{1}(\beta)N(\alpha)+sD(\alpha),t)ds,\label{E:cylinder}\\
y(\beta)  &  =(R\cos\beta,R\sin\beta).\nonumber
\end{align}
where we assume for simplicity that the detectors are infinitely long.

\begin{problem}
Reconstruct the initial condition $f(x)$ (supported within a ball of radius
$R$ centered at the origin) from the measurements $P_{\alpha}(y(\beta),t)$
known for all $\alpha\in\lbrack0,\pi],$ $\beta\in\lbrack0,2\pi],$ $t\in
\lbrack0,\infty).$
\end{problem}

We solve this problem in Section~\ref{S:linedet}.

\section{Fast algorithms for Problem 1\label{S:conventional}}

As mentioned in the Introduction, there exist a variety of the methods for the
solution of Problem 1 in 2D and 3D. However, none of the known methods have
optimal computational complexity. A fast algorithm for the 3D version of the
problem is needed since even relatively fast $\mathcal{O}(n^{4})$ methods
require several hours of computing time; methods based on slower
$\mathcal{O}(n^{5})$ discrete versions of explicit inversion formulas can
easily run a couple of days on larger computational grids.

A single solution of the 2D problem does not take long in practical terms even
if a slow algorithm is used. However, in order to solve Problem 2, one need to
solve the 2D version of Problem 1 $\mathcal{O}(n)$ times, which can again
result in several hours of computation. Hence, a fast method is also needed.

\subsection{2D case}

\label{S:2Dconventional}

\subsubsection{The algorithm\label{S:alg2d}}

Solution $u(x,t)$ to the initial value problem (\ref{IVP}) in 2D can be
written \cite{Vladimirov} in the following form:%
\[
u(y,t)=\int\limits_{B}f(x)\frac{\partial}{\partial t}\Phi(y-x,t)dx
\]
where $\Phi(x,t)$ is the free space Green function of the wave equation:
\[
\Phi(x,t)=\left\{
\begin{array}
[c]{cc}%
\frac{1}{2\pi\sqrt{t^{2}-x^{2}}}, & t>|x|,\\
0, & t\leq|x|.
\end{array}
\right.
\]
In particular, the measured data $P(y,t)$ are given by the similar formula
\begin{equation}
P(y,t)=\int\limits_{B}f(x)\frac{\partial}{\partial t}\Phi(y-x,t)dx,\qquad
y\in\partial B, \label{grconv}%
\end{equation}
where $B$ is disk of radius $R$ and $\partial B$ is its boundary. Let us find
the Fourier transform $\hat{P}(y,t)$ of $P(y,t)$ in $t:$%
\begin{equation}
\hat{P}(y,\lambda)\equiv\int\limits_{\mathbb{R}}P(y,t)e^{it\lambda
}dt=-i\lambda\int\limits_{B}f(x)\hat{\Phi}(y-x,\lambda),\qquad y\in\partial B,
\label{Phat2D}%
\end{equation}
and the Fourier transform $\hat{\Phi}(x,\lambda)$ of the Green function:
\[
\hat{\Phi}(x,\lambda)\equiv\int\limits_{\mathbb{R}}\Phi(x,t)e^{it\lambda
}dt=\frac{i}{4}H_{0}^{(1)}(\lambda|x|).
\]
By combining (\ref{grconv}) and (\ref{Phat2D}) we obtain:%
\begin{equation}
\lambda\int\limits_{B}f(x)H_{0}^{(1)}(\lambda|y-x|)dx=4\hat{P}(y,\lambda).
\label{Helmholtz}%
\end{equation}
Let us utilize the addition theorem for $H_{0}^{(1)}$ (see, for example
\cite{Colton}):%
\begin{align}
H_{0}^{(1)}(\lambda|y-x|)  &  =\sum\limits_{k=-\infty}^{\infty}H_{|k|}%
^{(1)}(\lambda R)J_{|k|}(\lambda r)e^{ik(\varphi-\theta)},\label{addition}\\
x  &  =r\hat{x},\quad\hat{x}(\theta)=(\cos\theta,\sin\theta),\nonumber\\
y(R,\varphi)  &  =R(\cos\varphi,\sin\varphi),\nonumber\\
R  &  >r.\nonumber
\end{align}
Also, let us expand $\hat{P}(y(R,\varphi),\lambda)$ and $f(r\hat{x}(\theta))$
in the Fourier series in $\varphi$ and $\theta$:%
\begin{align}
\hat{P}(y(R,\varphi),\lambda)  &  =\sum\limits_{k=-\infty}^{\infty}\hat{P}%
_{k}(\lambda)e^{ik\varphi},\nonumber\\
f(r\hat{x}(\theta))  &  =\sum\limits_{m=-\infty}^{\infty}f_{m}(r)e^{im\theta},
\label{fseries}%
\end{align}
where coefficients $\hat{P}_{k}(\lambda)$ and $f_{m}(r)$ are given by the
formulas%
\begin{align*}
\hat{P}_{k}(\lambda)  &  =\frac{1}{2\pi}\int\limits_{0}^{2\pi}\hat
{P}(z(R,\varphi),\lambda)e^{-ik\varphi}d\varphi,\\
f_{m}(r)  &  =\frac{1}{2\pi}\int\limits_{0}^{2\pi}f(r,\theta)e^{-im\varphi
}d\varphi.
\end{align*}

By substituting (\ref{addition}) into (\ref{Helmholtz}) and expanding in the
Fourier series one obtains%
\begin{align}
\hat{P}_{k}(\lambda)  &  =\frac{\lambda}{8\pi}\int\limits_{0}^{2\pi}\left[
\int\limits_{0}^{2\pi}\int\limits_{0}^{\infty}f(r\hat{x}(\theta))\sum
\limits_{k=-\infty}^{\infty}H_{|k|}^{(1)}(\lambda R)J_{|k|}(\lambda
r)e^{ik(\varphi-\theta)}rdrd\theta\right]  e^{-ik\varphi}d\varphi\nonumber\\
&  =\lambda H_{|k|}^{(1)}(\lambda R)\left(  \int\limits_{0}^{2\pi}%
\int\limits_{0}^{\infty}f(r\hat{x}(\theta))J_{|k|}(\lambda r)e^{-ik\theta
}rdrd\theta,\right)  \label{series}%
\end{align}
and further, by utilizing (\ref{fseries}):%
\[
\hat{P}_{k}(\lambda)=\frac{\pi}{2}\lambda H_{|k|}^{(1)}(\lambda R)\int
\limits_{0}^{\infty}f_{|k|}(r)J_{|k|}(\lambda r)rdr.
\]
The above formula relates $\hat{P}_{k}(\lambda)$ to the Hankel transform of
$f_{|k|}(r),$ and, since the latter transform is self-invertible,%
\begin{equation}
f_{k}(r)=\frac{2}{\pi}\int\limits_{0}^{\infty}\frac{\hat{P}_{k}(\lambda
)}{H_{|k|}^{(1)}(\lambda R)}J_{|k|}(\lambda r)d\lambda. \label{harmsol}%
\end{equation}
Formula (\ref{harmsol}) was first presented in \cite{AmbKuch}. A somewhat
similar expression was obtained in \cite{Norton1}; instead of $H_{|k|}%
^{(1)}(\lambda R)$ it contained the Bessel functions $J_{|k|}(\lambda R)$ in
the denominator, and a term that corresponds to the real part of our $\hat
{P}_{k}(\lambda)$ in the numerator. In that formula, theoretically, the zeros
of $\operatorname{Re}\hat{P}_{k}(\lambda)$ cancel the zeros of the Bessel
function in the denominator. However, when $\hat{P}_{k}(\lambda)$ is computed
from the data of real measurements contaminated by noise, the cancellation
will not happen automatically. In \cite{Haltm-circ} a technique based on the
use of Fourier-Bessel series was proposed to avoid such division by zeros. Our
approach is based on formula (\ref{harmsol}); since Hankel's functions do not
vanish for real (and bounded) values of the argument, this formula provides a
stable (and simple) way to recover $f_{k}(r).$

A straightforward discretization of equation (\ref{harmsol}) leads to a method
that would require $\mathcal{O}(n^{2})$ flops per each $f_{k}$ and
$\mathcal{O}(n^{3})$ per whole reconstruction, and thus, is not fast enough.
To develop a fast algorithm we combine (\ref{fseries}) and (\ref{harmsol}) and
represent $f$ in the following form:%
\begin{equation}
f(r\hat{x}(\theta))=\frac{2}{\pi}\sum\limits_{k=-\infty}^{\infty}\left(
\int\limits_{0}^{\infty}\frac{\hat{P}_{k}(\lambda)}{H_{|k|}^{(1)}(\lambda
R)}J_{|k|}(\lambda r)e^{ik\theta}d\lambda\right)  . \label{myseries}%
\end{equation}
Let us now consider plain waves $W(r\hat{x}(\theta),\Lambda)$:%
\begin{align*}
W(r\hat{x}(\theta),\Lambda)  &  =\exp(i\lambda r\cos(\theta-\varphi)),\\
\Lambda &  =\lambda(\cos\varphi,\sin\varphi),
\end{align*}
and expand these waves in the Fourier series in $\varphi.$ The values of the
corresponding Fourier coefficients can be found using the Jacobi-Anger
expansion \cite{Colton}:%
\[
W(r\hat{x}(\theta),\Lambda)=\sum\limits_{k=-\infty}^{\infty}i^{|k|}%
J_{|k|}(\lambda r)e^{ik\theta},
\]
so that%
\[
J_{|k|}(\lambda r)e^{ik\theta}=\frac{(-i)^{|k|}}{2\pi}\int\limits_{0}^{2\pi
}W(r\hat{x}(\theta),\Lambda(\lambda,\varphi))e^{ik\varphi}d\varphi.
\]
By substituting the latter formula in (\ref{myseries}) we obtain%
\begin{equation}
f(r\hat{x}(\theta))=\frac{1}{\pi^{2}}\int\limits_{0}^{\infty}\int
\limits_{0}^{2\pi}\left[  \sum\limits_{k=-\infty}^{\infty}\frac{(-i)^{k}%
\hat{P}_{k}(\lambda)}{\lambda H_{|k|}^{(1)}(\lambda R)}e^{ik\varphi}\right]
W(r\hat{x}(\theta),\Lambda(\lambda,\varphi))d\varphi\lambda d\lambda.
\label{long2d}%
\end{equation}
Let us denote by $\hat{f}(\Lambda(\lambda,\varphi))$ the expression in the
brackets (the choice of such a notation will become clear momentarily)%
\begin{equation}
\hat{f}(\Lambda(\lambda,\varphi))=\frac{2}{\pi}\sum\limits_{k=-\infty}%
^{\infty}\frac{(-i)^{k}\hat{P}_{k}(\lambda)}{\lambda H_{|k|}^{(1)}(\lambda
R)}e^{ik\varphi}. \label{fastalg}%
\end{equation}
Then (\ref{long2d}) can be re-written as%
\[
f(x)=\frac{1}{2\pi}\int\limits_{\mathbb{R}^{2}}\hat{f}(\Lambda)\exp
(ix\cdot\Lambda)d\Lambda,
\]
which means that $\hat{f}(\Lambda)$ is, in fact, the 2D Fourier transform of
$f(x)$ defined in the standard way:%
\[
\hat{f}(\Lambda)=\frac{1}{2\pi}\int\limits_{\mathbb{R}^{2}}f(x)\exp
(-ix\cdot\Lambda)dx.
\]

Formula (\ref{fastalg}) allows us to compute $\hat{f}(\Lambda)$ for all values
of $\Lambda\neq0.$ In order to find $\hat{f}(0)$ the following integral
identity can be used:%
\begin{align}
\hat{f}(0)  &  =\frac{1}{2\pi}\int\limits_{\mathbb{R}^{2}}f(x)dx=\int
\limits_{0}^{R}rf_{0}(r)dr=\int\limits_{0}^{R}r\left[  \int\limits_{0}%
^{\infty}\frac{2\hat{P}_{0}(\lambda)}{\pi H_{0}^{(1)}(\lambda R)}J_{0}(\lambda
r)d\lambda\right]  dr\nonumber\\
&  =\int\limits_{0}^{\infty}\frac{2\hat{P}_{0}(\lambda)}{\pi H_{0}%
^{(1)}(\lambda R)}\left[  \int\limits_{0}^{R}rJ_{0}(\lambda r)dr\right]
d\lambda=\int\limits_{0}^{\infty}\frac{2\hat{P}_{0}(\lambda)}{\pi\lambda
H_{0}^{(1)}(\lambda R)}RJ_{1}(\lambda R)d\lambda. \label{f0in2d}%
\end{align}
Since $H_{0}^{(1)}$ has a logarithmic singularity at $\lambda=0,$ and $J_{1}$
has a single root at the latter point, the integrand in (\ref{f0in2d})
vanishes as $\lambda\rightarrow0,$ and the formula is well-defined. Now $f(x)$
can be computed by applying the 2D inverse Fourier transform to $\hat
{f}(\Lambda).$ We utilize the FFT's on
various steps of the resulting algorithm to make the computations fast:

\textbf{The algorithm for Problem 1 in 2D:}

\begin{enumerate}
\item On an equispaced grid in $\lambda$ compute $\hat{P}(y,\lambda)$ using 1D
FFT in time:%
\[
\hat{P}(y,\lambda)=\int\limits_{\mathbb{R}}P(y,t)e^{it\lambda}dt.
\]

\item For each value of $\lambda$ in the grid compute $\hat{P}_{k}(\lambda)$
using 1D FFT in $\varphi$:
\[
\hat{P}_{k}(\lambda)=\frac{1}{2\pi}\int\limits_{0}^{2\pi}\hat{P}%
(y(R,\varphi),\lambda)e^{-ik\varphi}d\varphi.
\]

\item For each value of $\lambda>0$ in the grid compute coefficients
$b_{k}(\lambda)$%
\begin{equation}
b_{k}(\lambda)=\frac{2(-i)^{k}\hat{P}_{k}(\lambda)}{\pi\lambda H_{|k|}%
^{(1)}(\lambda R)}. \label{bkeys}%
\end{equation}

\item On a polar grid in $\lambda$ and $\varphi$ compute $\hat{f}%
(\Lambda(\lambda,\varphi))$ by summing the Fourier series (use 1D FFT):
\[
\hat{f}(\Lambda(\lambda,\varphi))=\sum\limits_{k=-\infty}^{\infty}%
b_{k}(\lambda)e^{ik\varphi},\lambda\neq0.
\]

\item Compute $\hat{f}(0)$ from formula (\ref{f0in2d}) using the trapezoid rule.

\item Interpolate $\hat{f}(\Lambda)$ to a Cartesian grid in $\Lambda.$

\item Reconstruct $f(x)$ by the 2D inverse FFT.
\end{enumerate}

It is not difficult to estimate the computational costs of the present method.
Assuming that the computational grid is of size $n\times n,$ and the number of
detectors and the size of the grids in $\lambda$ and $\varphi$ are
$\mathcal{O}(n),$ steps 1, 2, 4, and 7 require $\mathcal{O}(n^{2}\log n)$
flops each. Step 3 needs $\mathcal{O}(n^{2})$ operations and step 5 is
completed in $\mathcal{O}(n)$ flops.

\begin{figure}[t]
\begin{center}
\subfigure[]{\includegraphics[width=1.8in,height=1.8in]{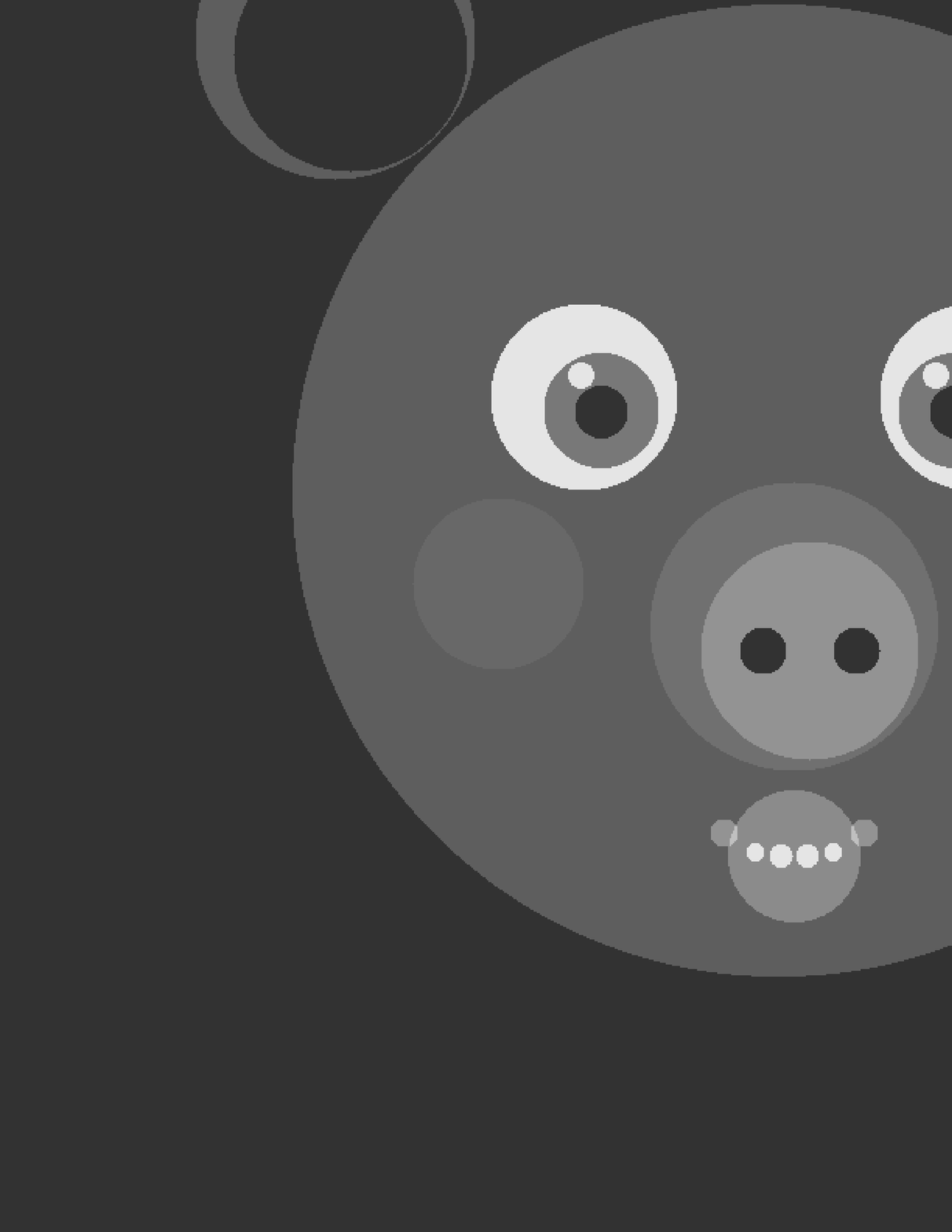} }
\subfigure[]{\includegraphics[width=1.8in,height=1.8in]{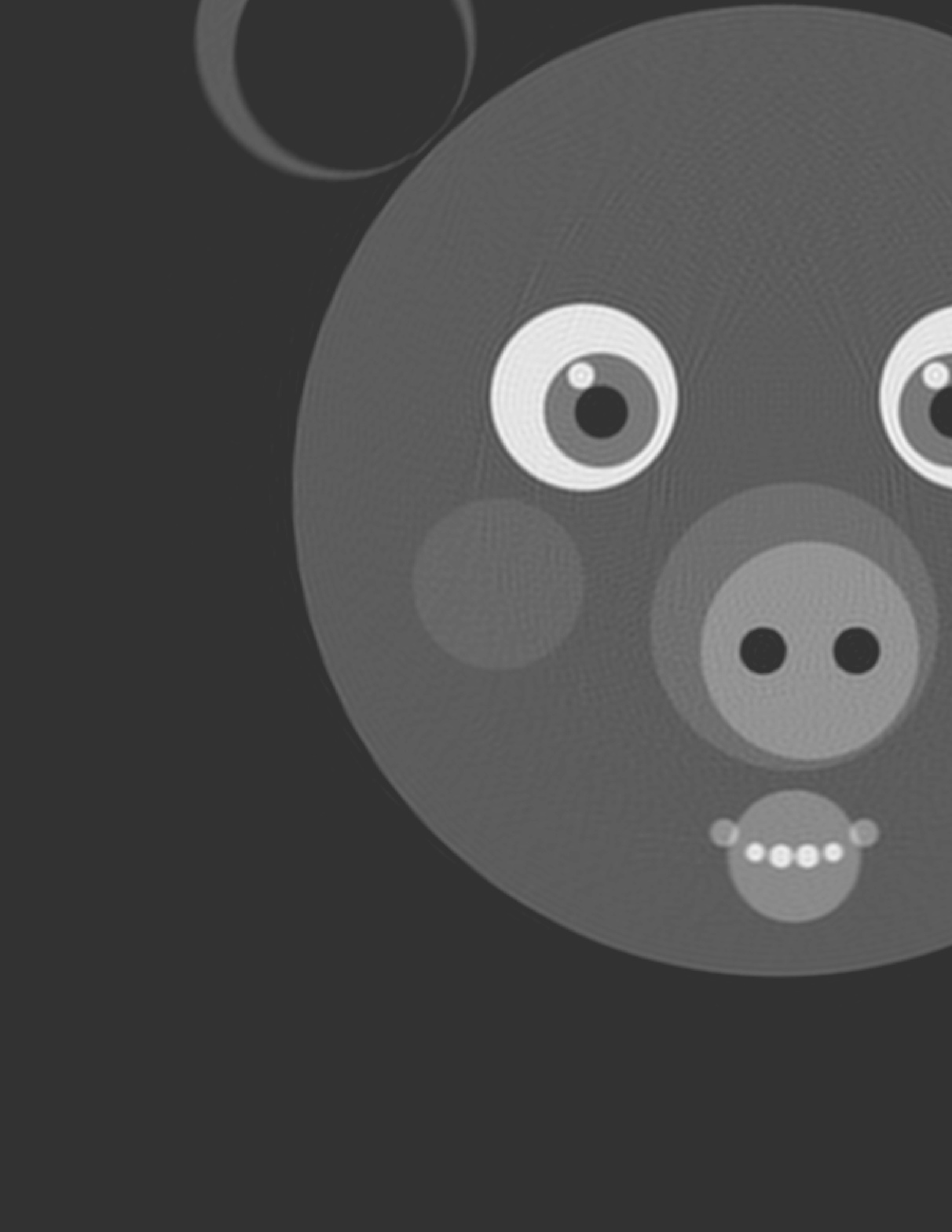} }
\subfigure[]{\includegraphics[width=1.8in,height=1.8in]{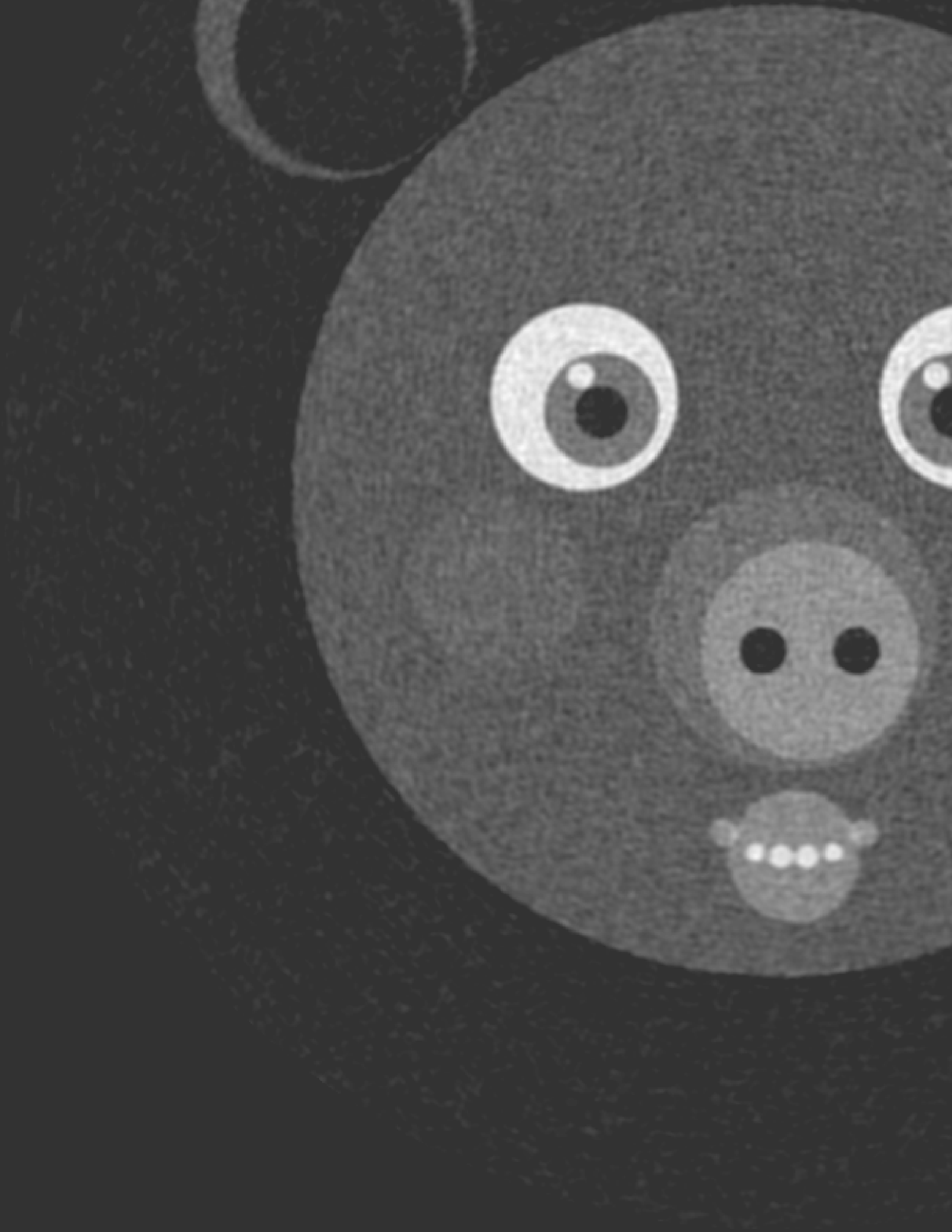}}
\newline%
\subfigure[]{\includegraphics[width=2.7in,height=1.6in]{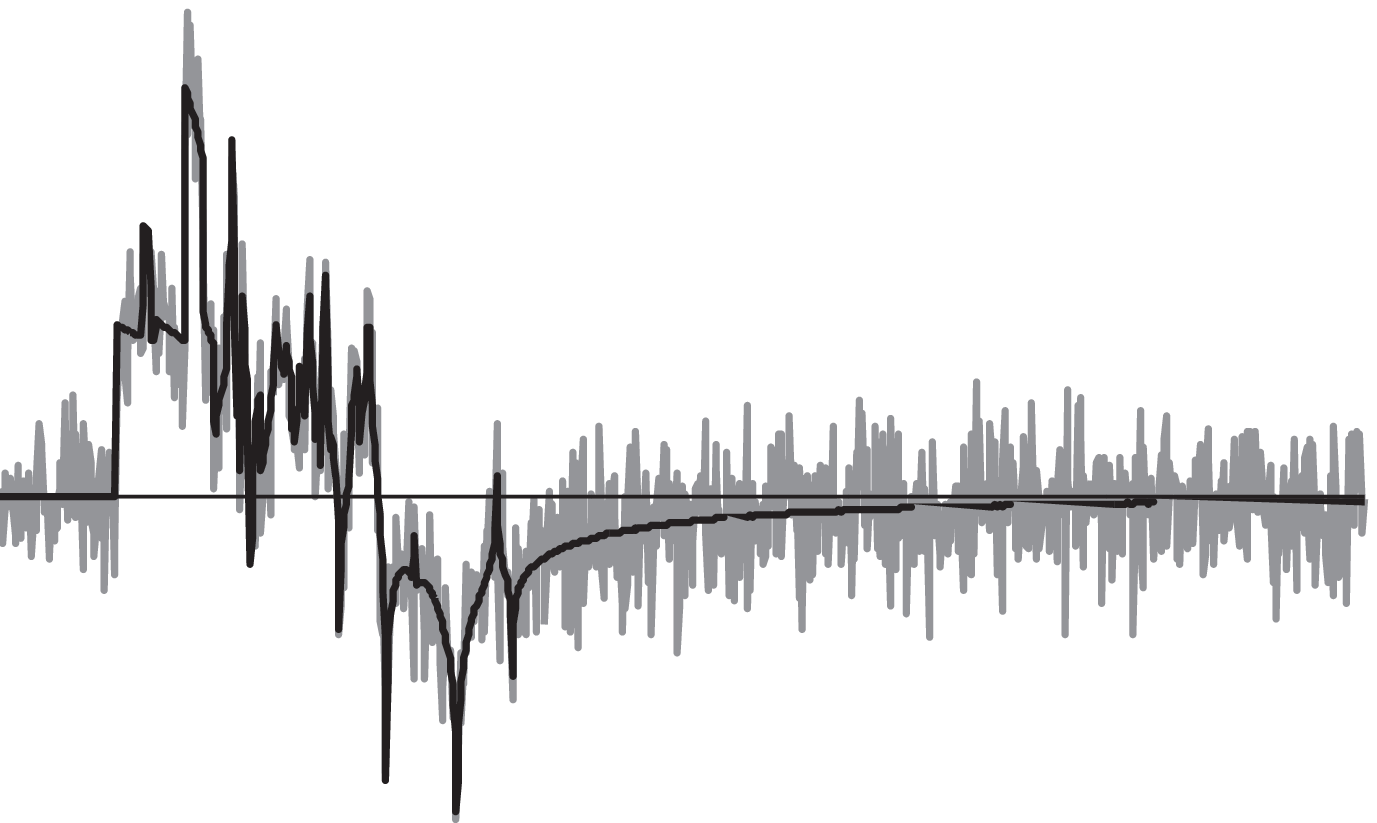}}
\end{center}
\caption{Reconstruction in $2D$ (a) phantom (b) reconstruction from noiseless
data (c) reconstruction in the presence of $50\%$ noise (d) comparison of
noiseless (black line) and noisy (gray line) data }%
\label{F:rec2D}%
\end{figure}

The interpolation step 6 needs some commentary. It is known \cite{Natterer}
that low-order interpolation in spectral domain can be a source of significant
error. In particular, interpolation in the radial direction (i.e. in
$\lambda)$ is more sensitive than that in the angular direction (i.e. in
$\varphi)$. The results presented below were obtained by combining linear
interpolation in $\varphi$ with interpolation by cubic polynomials in
$\lambda.$ We found that this yields accuracy that is more than sufficient for
practical tomography applications. Since such interpolation is local, the
associated operation count is $\mathcal{O}(n^{2}).$ However, if a more
accurate interpolation is desired, one can apply methods based on global
trigonometric interpolation using the NUFFT (see, e.g. \cite{Dutt}), for the
total cost of $\mathcal{O}(n^{2}\log n)$ operations. In any case, the
asymptotic estimate for the total number of operations required by the present
method is $\mathcal{O}(n^{2}\log n).$

\subsubsection{Numerical simulations}

In order to evaluate the performance of our algorithm we simulated high
resolution projections that corresponded to a phantom consisting of several
characteristic functions of disks lying inside a unit circle, as shown in
Figure~\ref{F:rec2D}(a). The number of the simulated detectors was equal to
272 (this corresponds to the number of detectors in the set of real
measurements presented in Section~\ref{S:real}). The detectors were placed on
a circle of radius 1.05. The measurements were simulated for the time interval
$[0,5]$, with time step 0.005 (i.e., 1000 time samples were simulated). A
smooth cut-off was applied at the end of this time interval, and the signal
for $t>5$ was neglected.

The result of the reconstruction on a grid $1000\times1000$ is shown in
Figure~\ref{F:rec2D}(b); the computation took 0.3 seconds on a desktop
computer with a $2.4$ GHz Intel Core 2 Duo processor. The code was written in
Fortran-95, and computations were not parallelized.

\begin{figure}[t]
\begin{center}
\subfigure[]{\includegraphics[width=1.8in,height=1.8in]{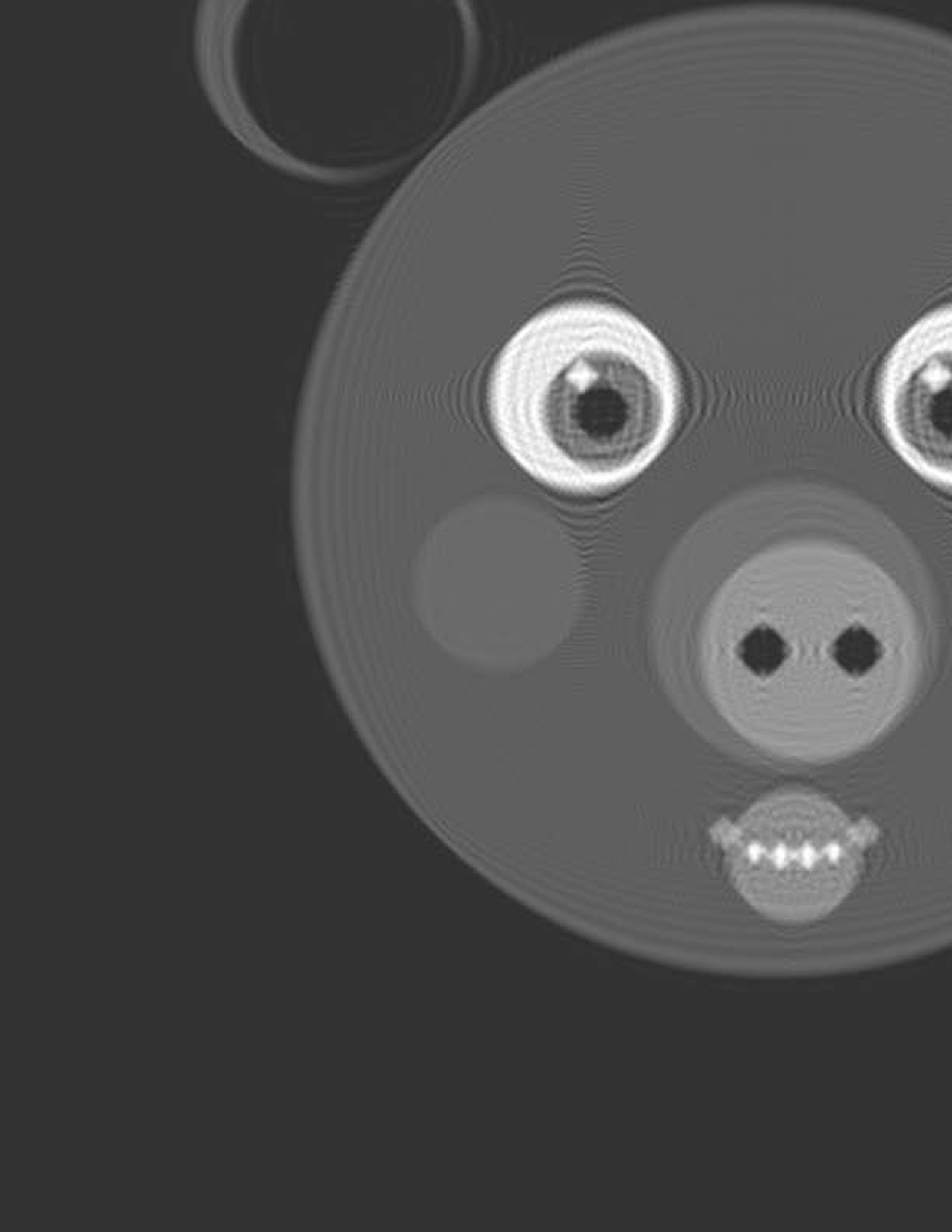}}
\subfigure[]{\includegraphics[width=1.8in,height=1.8in]{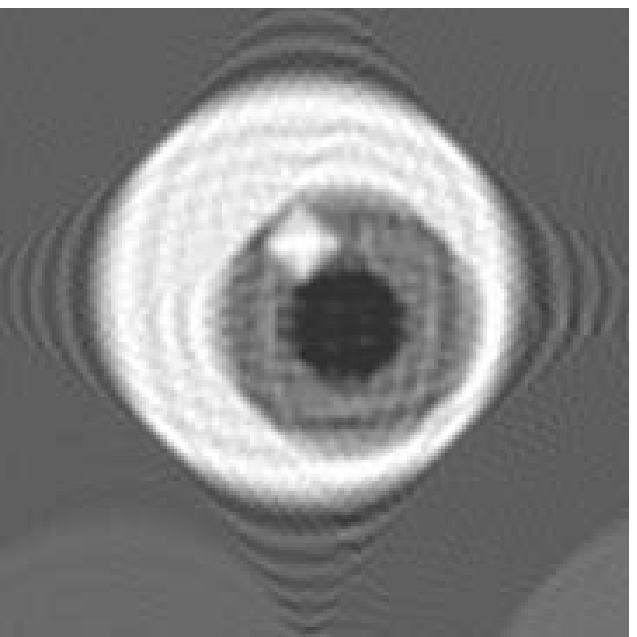}}
\subfigure[]{\includegraphics[width=1.8in,height=1.8in]{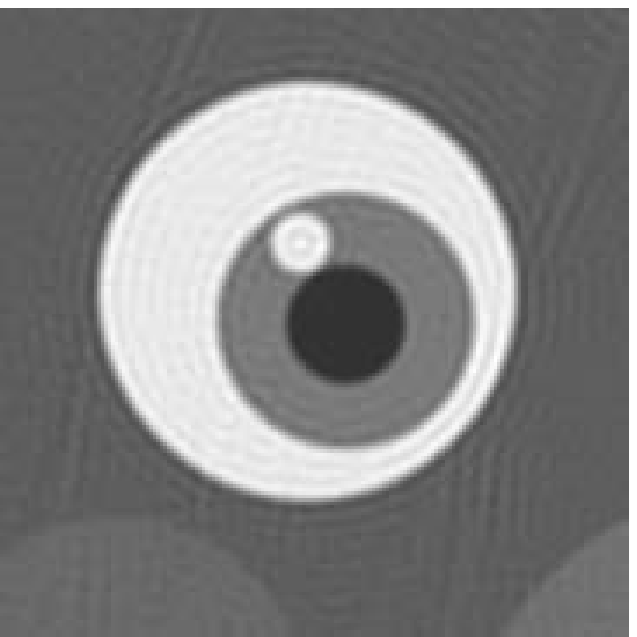}}
\end{center}
\caption{Comparison with the time reversal using finite differences: (a) image
obtained by the time reversal; (b) a fragment of the image in part (a); (c) a
fragment of the image in Figure~\ref{F:rec2D}(b) reconstructed by the present
method }%
\label{F:reversal}%
\end{figure}

In order to test stability of the present method we added $50\%$ (in $L^{2}$
norm) noise to the simulated projections. Figure~\ref{F:rec2D}(c) compares the
noiseless and noisy data for the first detector. Figure~\ref{F:rec2D}(d)
presents the result of the reconstruction from the noisy data. One can notice
low level of noise in the reconstructed image, which is unusual for the
inverse problems of tomography, where at least mild instability (and hence,
some noise amplification) occurs almost always. A quick glance at the steps of
the algorithm confirms that all the operations are stable: the Fourier
transforms and series are isometries in $L^{2}$, and in the
formula~(\ref{bkeys}) $\hat{P}_{k}(\lambda)$ is multiplied by a factor that
vanishes as $\lambda$ goes to infinity. The reason of such stability is the
presence of the time derivative in the equation (\ref{grconv}) that describes
the wave propagation. In the absence of this derivative (i.e. for example if
the initial value of the pressure was zero but the time derivative was
non-zero) the stability properties of the problem would be similar to those of
the classical inverse Radon transform in 2D, and some noise amplification
would occur during the reconstruction.

It is also interesting to compare the present technique with the results
obtained by the time reversal computed using finite differences. In order to
do so, we solved the 2D wave equation back in time in the unit square $\Omega
$, using the explicit leapfrog scheme%
\[
\frac{u_{k,l}^{m+1}-2u_{k,l}^{m}+u_{k,l}^{m-1}}{(\Delta t)^{2}}=\frac
{u_{k+1,l}^{m}+u_{k-1,l}^{m}+u_{k,l+1}^{m}+u_{k,l-1}^{m}-4u_{k,l}^{m}}{(\Delta
x)^{2}},
\]
starting with $u=0,\frac{\partial}{\partial t}u=0$ at time $t=5,$ and
enforcing the boundary condition $u(y,t)=P(y,t)$ for all $y$ lying on the
boundary $\partial\Omega$ of $\Omega.$ We chose to do the reconstruction in a
square domain since the exact boundary conditions are easy to enforce in the
nodes of the computational grids lying on the sides of the square. (For a
domain of a different shape enforcing the boundary conditions requires
interpolation (see e.g. \cite{burg-exac-appro}) which brings additional
error). In order to guarantee the stability of this explicit scheme we had to
increase the number of time steps from 1000 to 3800.

The result of the reconstruction by the time reversal is shown in
Figure~\ref{F:reversal}(a). Larger features of the image are reconstructed
quite well. However, the second-order finite difference scheme we used is not
very accurate on higher spatial frequencies, which leads to a typical
diamond-shaped distortion of small round shapes, clearly visible in the
magnified fragment of the image presented in Figure~\ref{F:reversal}(b). The
same part of the image in Figure~\ref{F:rec2D}(b) (obtained using the present
algorithm) is shown in Figure~\ref{F:reversal}(c). The artifacts are much
smaller, in spite of the lesser number of time steps in the data.

The time reversal took 214 seconds on the same computer, as before. Even if
one discounts this number by a factor of 4 (to account for a larger number of
the time steps), our technique is still faster by two orders of magnitude.

\subsection{3D case}

In this section we briefly outline the 3D version of the method presented in
Section~\ref{S:alg2d}. In 3D, the Green's function for the wave equation
satisfying radiation condition at infinity has the following form
\cite{Vladimirov}:%
\[
\Phi(x,t)=\frac{\delta(t-|x|)}{4\pi|x|}.
\]
As in 2D case, we compute the Fourier transform in $t$ from the data $P(y,t)$
\begin{equation}
\hat{P}(y,\lambda)\equiv\int\limits_{\mathbb{R}}P(y,t)e^{it\lambda
}dt=-i\lambda\int\limits_{B}f(x)\hat{\Phi}(y-x,\lambda),\qquad y\in\partial B,
\label{grconv3d}%
\end{equation}
and the Fourier transform $\hat{\Phi}(x,\lambda)$ from the Green's function%
\begin{equation}
\hat{\Phi}(x,\lambda)=\int\limits_{\mathbb{R}}\frac{\delta(t-|x|)}{4\pi
|x|}e^{it\lambda}dt=\frac{e^{i\lambda|x|}}{4\pi|x|}=\frac{i\lambda}{4\pi}%
h_{0}^{(1)}(\lambda|x|). \label{grfourier3d}%
\end{equation}
By combining (\ref{grconv3d}) and (\ref{grfourier3d})\ we obtain:%
\begin{equation}
\lambda^{2}\int\limits_{\Omega}f(x)h_{0}^{(1)}(\lambda|y-x|)dx=4\pi\hat
{P}(y,\lambda). \label{newconv3d}%
\end{equation}
We will utilize spherical harmonics $Y_{k}^{m}(\hat{z}),$ $\hat{z}%
\in\mathbb{S}^{2},$ normalized so that%
\[
\int\limits_{\mathbb{S}^{2}}Y_{k}^{m}(\hat{z})\overline{Y_{p}^{s}(\hat{z}%
)}d\hat{z}=\delta_{k,p}\delta_{m,s}.
\]
where $\delta_{k,p}$ is the Kronecker symbol. Let us extend $f(x)$ and
$\hat{P}(y,\lambda)$ in spherical harmonics:%

\begin{align}
f(r\hat{x})  &  =\sum\limits_{s=0}^{\infty}\sum\limits_{p=-s}^{s}%
f_{s,p}(r)Y_{s}^{p}(\hat{x}),\hat{x}\in\mathbb{S}^{2},r=|x|,\label{sphrhar1}\\
f_{s,p}(r)  &  =\int\limits_{\mathbb{S}^{2}}f(\hat{x}r)\overline{Y_{s}%
^{p}(\hat{x})}d\hat{x},\nonumber\\
\hat{P}(R\hat{y},\lambda)  &  =\sum\limits_{s=0}^{\infty}\sum\limits_{p=-s}%
^{s}\hat{P}_{s,p}(\lambda)Y_{s}^{p}(\hat{y}),\label{sphrhar3}\\
\hat{P}_{s,p}(\lambda)  &  =\int\limits_{\mathbb{S}^{2}}\hat{P}(R\hat
{y},\lambda)\overline{Y_{s}^{p}(\hat{y})}d\hat{y}. \label{sphrhar4}%
\end{align}

As in the 2D case, we make use of the addition theorem for $h_{0}^{(1)}$
\cite{Colton}:%
\[
h_{0}^{(1)}(\lambda|z-x|)=4\pi\sum\limits_{k=0}^{\infty}\sum\limits_{m=-k}%
^{k}h_{k}^{(1)}(\lambda|z|)j_{k}(\lambda|x|)Y_{k}^{m}(\hat{z})\overline
{Y_{k}^{m}(\hat{x})},
\]
where $j_{k}(\cdot)$ and $h_{k}^{(1)}(\cdot)$ are, respectively, the spherical
Bessel and Hankel functions.

By substituting the above equation, together with (\ref{sphrhar1}), into
(\ref{newconv3d}) we obtain%
\begin{align*}
\hat{P}(y,\lambda)  &  =\int\limits_{0}^{\infty}\int\limits_{\mathbb{S}^{2}%
}\left(  \sum\limits_{s=0}^{\infty}\sum\limits_{p=-s}^{s}f_{s,p}(r)Y_{s}%
^{p}(\hat{x})\right)  \left[  \sum\limits_{k=0}^{\infty}\sum\limits_{m=-k}%
^{k}h_{k}^{(1)}(\lambda R)j_{k}(\lambda r)Y_{k}^{m}(\hat{y})\overline
{Y_{k}^{m}(\hat{x})}\right]  d\hat{x}r^{2}dr\\
&  =\sum\limits_{k=0}^{\infty}\sum\limits_{m=-k}^{k}\left(  h_{k}%
^{(1)}(\lambda R)\lambda^{2}\int\limits_{0}^{\infty}f_{k,m}(r)j_{k}(\lambda
r)r^{2}dr\right)  Y_{k}^{m}(\hat{y}),
\end{align*}
and further, by comparing with (\ref{sphrhar3}):%
\begin{equation}
\hat{P}_{s,p}(\lambda)=\lambda^{2}h_{s}^{(1)}(\lambda R)\left(  \int
\limits_{0}^{\infty}f_{s,p}(r)j_{s}(\lambda r)r^{2}dr\right)  .
\label{sphrtran}%
\end{equation}
The spherical Bessel functions are related to their cylindrical counterparts
by the equation%
\[
j_{s}(t)=\sqrt{\frac{\pi}{2t}}J_{s+1/2}(t).
\]
Thus, $\hat{P}_{s,p}(\lambda)$ can be expressed in terms of the Hankel
transform of $\sqrt{r}f_{s,p}(r)$:
\[
\hat{P}_{s,p}(\lambda)=\sqrt{\frac{\pi}{2}}\lambda^{3/2}h_{s}^{(1)}(\lambda
R)\left(  \int\limits_{0}^{\infty}\left[  \sqrt{r}f_{s,p}(r)\right]
J_{s+1/2}(\lambda r)rdr\right)  .
\]
Since the Hankel transforms are self-invertible, one recovers $f_{s,p}(r)$ as
follows:%
\begin{equation}
f_{s,p}(r)=\sqrt{\frac{2}{\pi}}\int\limits_{0}^{\infty}\frac{\hat{P}%
_{s,p}(R,\lambda)}{\lambda^{3/2}r^{1/2}h_{s}^{(1)}(\lambda R)}J_{s+1/2}%
(\lambda r)\lambda d\lambda=\frac{2}{\pi}\int\limits_{0}^{\infty}\frac{\hat
{P}_{s,p}(R,\lambda)}{h_{s}^{(1)}(\lambda R)}j_{s}(\lambda r)d\lambda.
\label{nortonsol3d}%
\end{equation}

This approach is close (although not quite identical) to the solution obtained
in \cite{Norton2}. (An expression equivalent to (\ref{nortonsol3d}) was
derived in \cite{Ramm}.) A straightforward computation of (\ref{nortonsol3d})
for all $s\leq n,$ $|p|\leq s$ leads to an algorithm of complexity at least
$\mathcal{O}(n^{4})$ for an $n\times n\times n$ computational grid. In order
to accelerate the computations we choose the following approach. By
substituting (\ref{nortonsol3d}) into (\ref{sphrhar1}) we obtain:%
\begin{equation}
f(r\hat{x})=\frac{2}{\pi}\int\limits_{0}^{\infty}\sum\limits_{s=0}^{\infty
}\sum\limits_{p=-s}^{s}\frac{\hat{P}_{s,p}(R,\lambda)}{h_{s}^{(1)}(\lambda
R)}j_{s}(\lambda r)Y_{s}^{p}(\hat{x})d\lambda. \label{myser3d}%
\end{equation}
A convenient integral representation for the term $j_{s}(\lambda r)Y_{s}%
^{p}(\hat{x})$ in the above equation is given by the Funk-Hecke formula
\cite{Colton}:
\begin{equation}
j_{s}(\lambda r)Y_{s}^{p}(\hat{x})=\frac{i^{s}}{4\pi}\int\limits_{\mathbb{S}%
^{2}}e^{-i\lambda x\cdot\hat{z}}Y_{s}^{p}(\hat{z})d\hat{z}. \label{Funk3d}%
\end{equation}
Combining (\ref{myser3d})\ and (\ref{Funk3d}) yields:%
\begin{align*}
f(x)  &  =\frac{1}{2\pi^{2}}\int\limits_{0}^{\infty}\int\limits_{\mathbb{S}%
^{2}}\sum\limits_{s=0}^{\infty}\sum\limits_{p=-s}^{s}\frac{i^{s}\hat{P}%
_{s,p}(R,\lambda)}{\lambda^{2}h_{s}^{(1)}(\lambda R)}e^{-i\lambda x\cdot
\hat{z}}Y_{s}^{p}(\hat{z})d\hat{z}\lambda^{2}d\lambda\\
&  =\frac{1}{(2\pi)^{3/2}}\int\limits_{\mathbb{R}^{3}}\left[  \sum
\limits_{s=0}^{\infty}\sum\limits_{p=-s}^{s}\sqrt{\frac{2}{\pi}}\frac
{i^{s}\hat{P}_{s,p}(R,|\Lambda|)}{\Lambda^{2}h_{s}^{(1)}(|\Lambda|R)}Y_{s}%
^{p}\left(  \frac{\Lambda}{|\Lambda|}\right)  \right]  e^{-ix\cdot\Lambda
}d\Lambda,
\end{align*}
where $\Lambda=\hat{z}\lambda.$ Let us denote the term in the brackets by
$F(\Lambda)$:%
\begin{equation}
F(\Lambda)=\sum\limits_{s=0}^{\infty}\sum\limits_{p=-s}^{s}\sqrt{\frac{2}{\pi
}}\frac{i^{s}\hat{P}_{s,p}(R,|\Lambda|)}{\Lambda^{2}h_{s}^{(1)}(|\Lambda
|R)}Y_{s}^{p}\left(  \frac{\Lambda}{|\Lambda|}\right)  . \label{fourierin3d}%
\end{equation}
Then%
\begin{equation}
f(x)=\frac{1}{(2\pi)^{3/2}}\int\limits_{\mathbb{R}^{3}}F(\Lambda
)e^{-ix\cdot\Lambda}d\Lambda, \label{recon3d}%
\end{equation}
and it becomes clear that function $F(\Lambda)$ is the 3D Fourier transform of
$f(x)$ defined as follows:%
\[
F(\Lambda)=\frac{1}{(2\pi)^{3/2}}\int\limits_{\mathbb{R}^{3}}f(x)e^{ix\cdot
\Lambda}dx.
\]
\bigskip

Formulas (\ref{grconv3d}) and (\ref{sphrhar4}) allow us to reconstruct the
Fourier transform $F(\Lambda)$ for all $\Lambda\neq0.$ In order to find $F(0)$
we notice that%
\[
F(0)=\frac{1}{(2\pi)^{3/2}}\int\limits_{\mathbb{R}^{3}}f(x)dx=\frac{1}%
{(2\pi)^{3/2}}\int\limits_{0}^{R}f_{0,0}(r)r^{2}dr,
\]
where $f_{0,0}$ is given by (\ref{nortonsol3d}), so that:%
\begin{align*}
F(0)  &  =\frac{1}{\sqrt{2}(\pi)^{5/2}}\int\limits_{0}^{R}\frac{\hat{P}%
_{0,0}(R,\lambda)}{h_{0}^{(1)}(\lambda R)}j_{0}(\lambda r)d\lambda r^{2}dr\\
&  =\frac{1}{\sqrt{2}(\pi)^{5/2}}\int\limits_{0}^{\infty}\frac{\hat{P}%
_{0,0}(R,\lambda)}{h_{0}^{(1)}(\lambda R)}\left[  \int\limits_{0}^{R}%
j_{0}(\lambda r)r^{2}dr\right]  d\lambda.
\end{align*}
Since $j_{0}(t)=\sin(t)/t,$ the integral in the brackets can be easily
evaluated:%
\[
\int\limits_{0}^{R}j_{0}(\lambda r)r^{2}dr=\frac{R}{\lambda^{2}}\left(
\frac{\sin\lambda R}{\lambda R}-\cos\lambda R\right)  ,
\]
leading to the following expression for $F(0)$:%
\begin{align}
F(0)  &  =\frac{R}{\sqrt{2}(\pi)^{5/2}}\int\limits_{0}^{\infty}\frac{\hat
{P}_{0,0}(R,\lambda)}{\lambda^{2}h_{0}^{(1)}(\lambda R)}\left(  \frac
{\sin\lambda R}{\lambda R}-\cos\lambda R\right)  d\lambda\nonumber\\
&  =\frac{iR^{2}}{\sqrt{2}(\pi)^{5/2}}\int\limits_{0}^{\infty}\frac{1}%
{\lambda}\hat{P}_{0,0}(R,\lambda)\exp(-i\lambda R)\left(  \frac{\sin\lambda
R}{\lambda R}-\cos\lambda R\right)  d\lambda. \label{f03d}%
\end{align}
Thus, $f(x)$ can be reconstructed by the following method:

\textbf{The algorithm for Problem 1 in 3D:}

\begin{enumerate}
\item On an equispaced grid in $\lambda$ compute $\hat{P}(y,\lambda)$ using 1D
FFT in time:%
\[
\hat{P}(y,\lambda)=\int\limits_{\mathbb{R}}P(y,t)e^{it\lambda}dt.
\]

\item For each value of $\lambda$ in the grid expand $\hat{P}(R\hat{y}%
,\lambda)$ in spherical harmonics in $\hat{y}$:
\[
\hat{P}_{s,p}(\lambda)=\int\limits_{\mathbb{S}^{2}}\hat{P}(R\hat{y}%
,\lambda)\overline{Y_{s}^{p}(\hat{y})}d\hat{y}.
\]

\item For each value of $\lambda>0$ in the grid compute coefficients
$b_{s,p}(\lambda)$:%
\[
b_{s,p}(\lambda)=\sqrt{\frac{2}{\pi}}\frac{i^{s}\hat{P}_{s,p}(R,\lambda
)}{\lambda^{2}h_{s}^{(1)}(\lambda R)}.
\]

\item On a spherical grid in $\Lambda$ compute $F(\Lambda)$ by summing
spherical harmonics:
\[
F(\Lambda)=\sum\limits_{s=0}^{\infty}\sum\limits_{p=-s}^{s}b_{s,p}%
(|\Lambda|)Y_{s}^{p}\left(  \frac{\Lambda}{|\Lambda|}\right)  ,\quad
\Lambda\neq0.
\]

\item Compute $F(0)$ from formula (\ref{f03d}) using the trapezoid rule.

\item Interpolate $F(\Lambda)$ to a Cartesian grid in $\Lambda.$

\item Reconstruct $f(x)$ by the 3D inverse FFT.
\end{enumerate}

The most time-consuming steps of this algorithm are steps 2 and 4,
corresponding to the Fourier analysis and synthesis on a sphere. One of the
simplest ways to decompose a function in the spherical harmonics combines
expansion into the Fourier series (which can be done fast using the FFT) with
the expansion into the series of Legendre functions. The straightforward
implementation of the latter Legendre transform requires $\mathcal{O}(n^{2})$
flops for a one dimensional function defined by $n$ samples. Then, expanding
one 2D function on a sphere (defined by $n\times n$ samples ) into a series of
spherical harmonics needs $\mathcal{O}(n^{3})$ operations, and the present
method involving $\mathcal{O}(n)$ of such expansions would require
$\mathcal{O}(n^{4})$ flops (synthesis of spherical harmonics is quite similar
to the analysis).
\begin{figure}[t]
\begin{center}
\includegraphics[width=2.8in,height=2.6in]{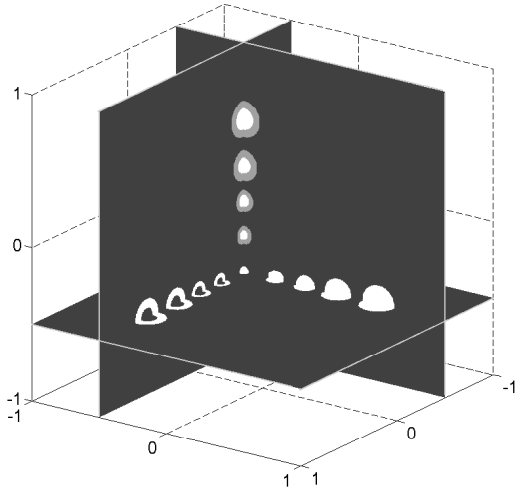}
\end{center}
\caption{3D phantom}%
\label{F:3Dphan}%
\end{figure}

Recently, several versions of the Fast Spherical Harmonics (FSH) transform
were introduced by several research groups \cite{Healy,Suda,Mohlen,Potts}. FSH
is asymptotically fast; it performs Fourier analysis and synthesis on a sphere
in $\mathcal{O}(n^{2}\log^{2}n)$ flops. If this algorithm is used on steps 2
and 4 of our reconstruction technique, the resulting method will also be fast,
requiring $\mathcal{O}(n^{3}\log^{2}n)$ operations per a 3D reconstruction.
One has to keep in mind, however, that the break-even size for which the
FSH\ starts to outperform the simple slow method described in the previous
paragraph is currently of order of hundreds, and thus our 3D reconstruction
method would not significantly outperform slower $\mathcal{O}(n^{4})$
techniques for the sizes of the problems we consider in this paper. In
addition, ready-to-use implementations of the FSH are only available for a
restricted number of operating systems and computational languages. However,
as the work on the FSH progresses, the performance and the ease of programming
of our 3D method will improve.

\section{Fast algorithm for Problem 2\label{S:linedet}}

It is known \cite{Palt-Machzend,burg-FBP} that Problem 2 can be reduced to
solving a set of Problems 1 in 2D, followed by a set of numerical inversions
of the 2D Radon transform. Indeed, let us fix unit vector $D(\alpha)$ and
consider all line integrals $v_{\alpha}(h,t)$ in the direction $D(\alpha)$ of
$u(x,t)$ :%
\begin{align*}
v_{\alpha}(h,t)  &  =\int\limits_{\mathbb{R}^{1}}u(h_{1}N(\alpha)+h_{2}%
e_{2}+sD(\alpha),t)ds,\\
h  &  =(h_{1},h_{2}).
\end{align*}
\begin{figure}[h]
\begin{center}
\subfigure{\includegraphics[width=1.8in,height=1.8in]{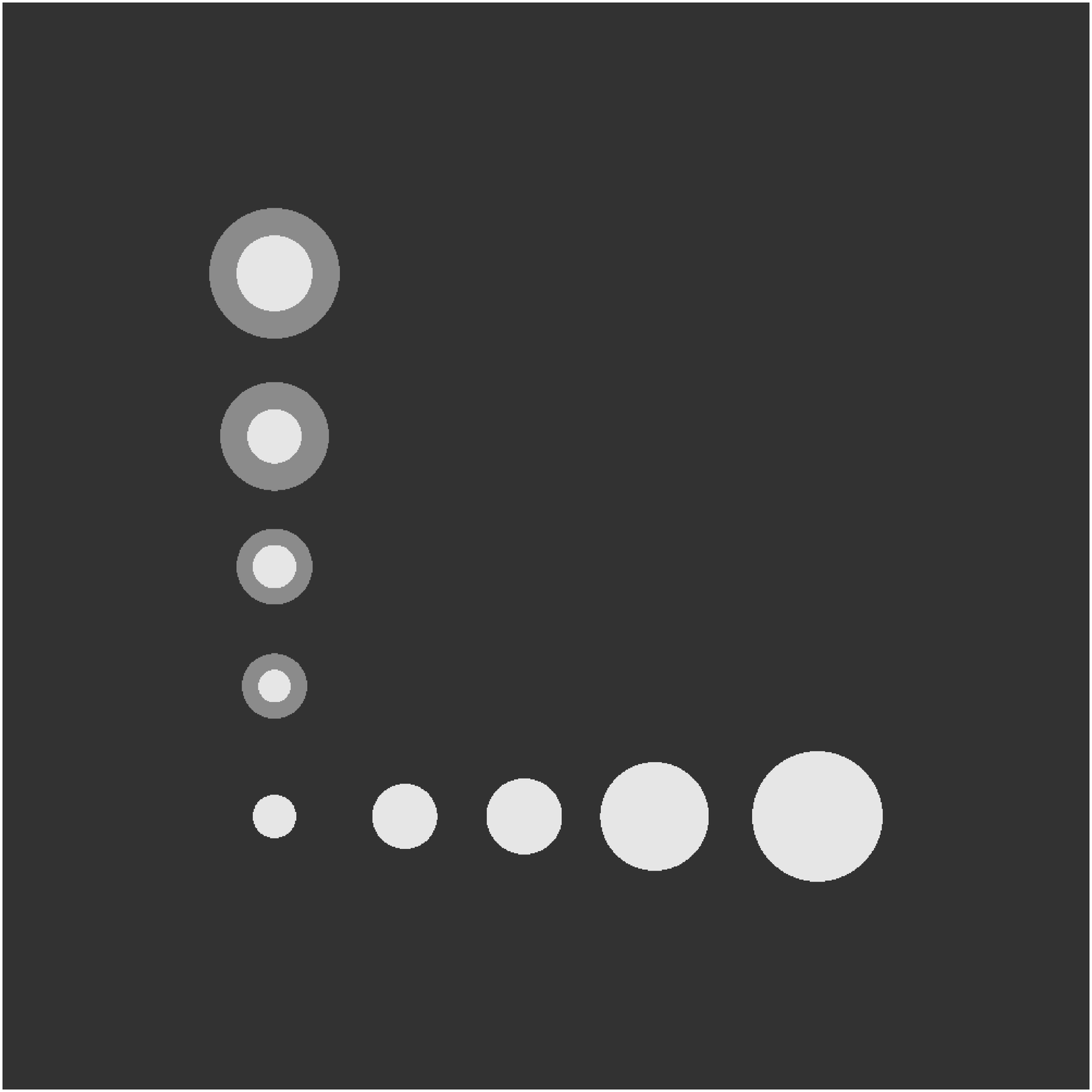}}
\subfigure{\includegraphics[width=1.8in,height=1.8in]{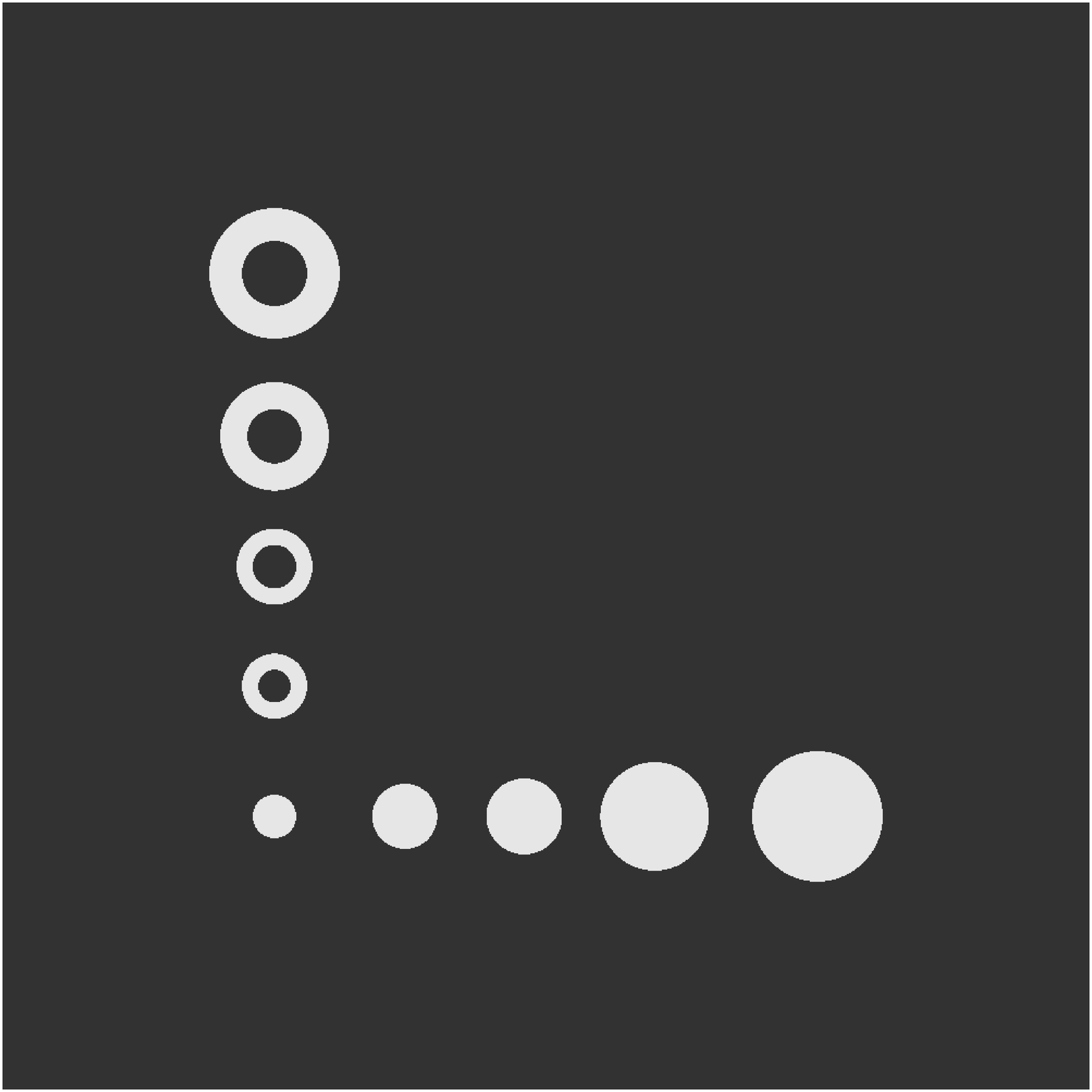}}
\subfigure{\includegraphics[width=1.8in,height=1.8in]{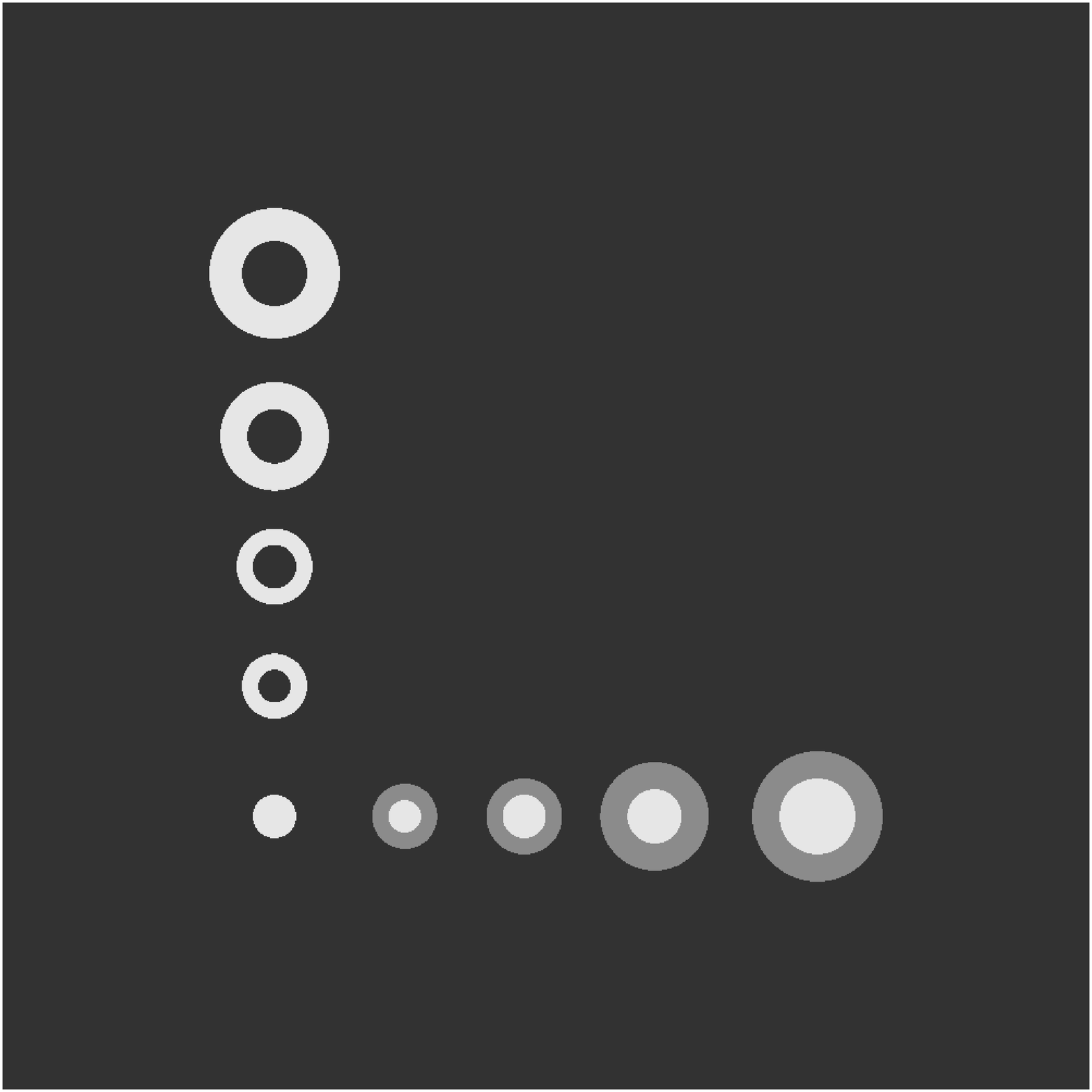}}\newline%
\subfigure{\includegraphics[width=1.8in,height=1.8in]{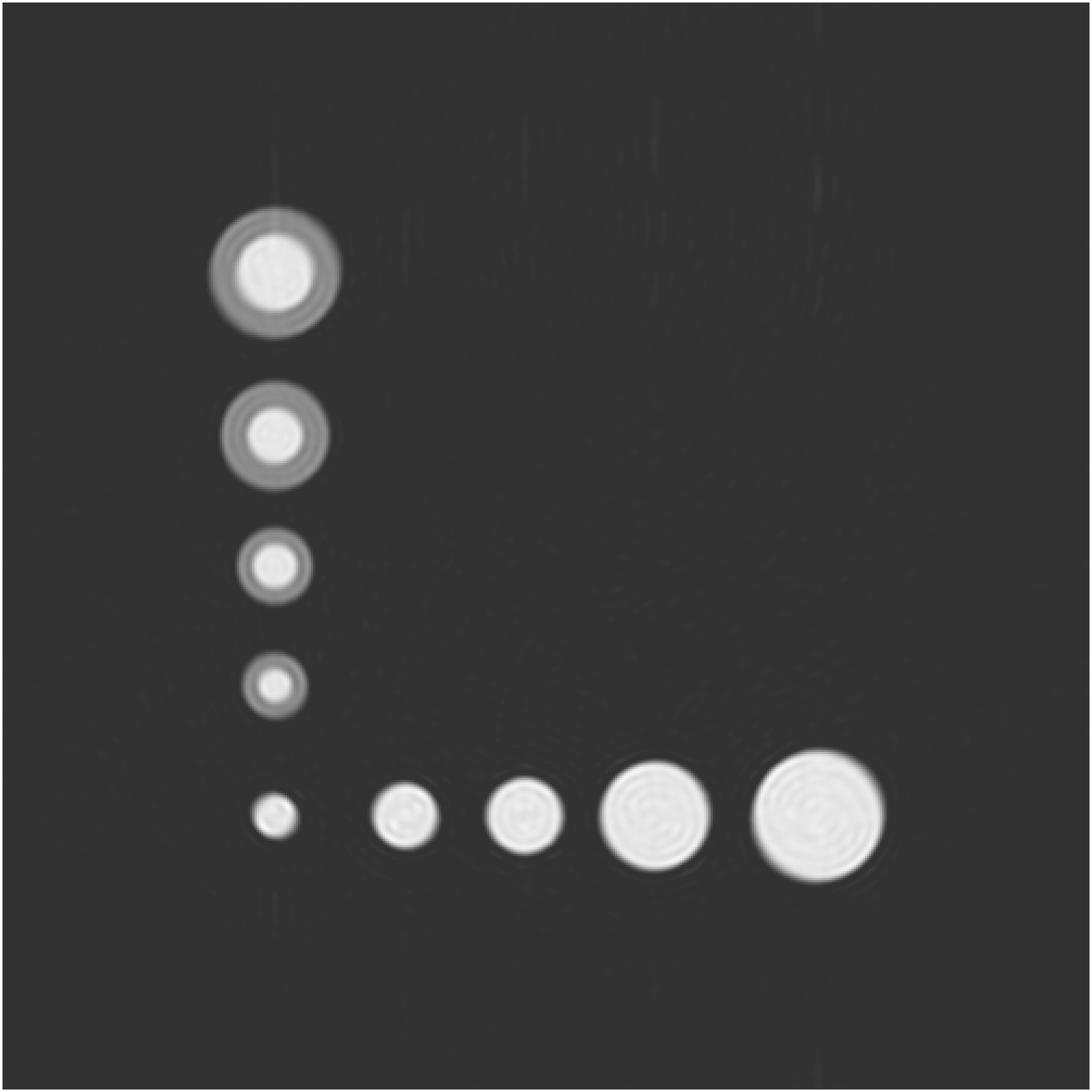}}
\subfigure{\includegraphics[width=1.8in,height=1.8in]{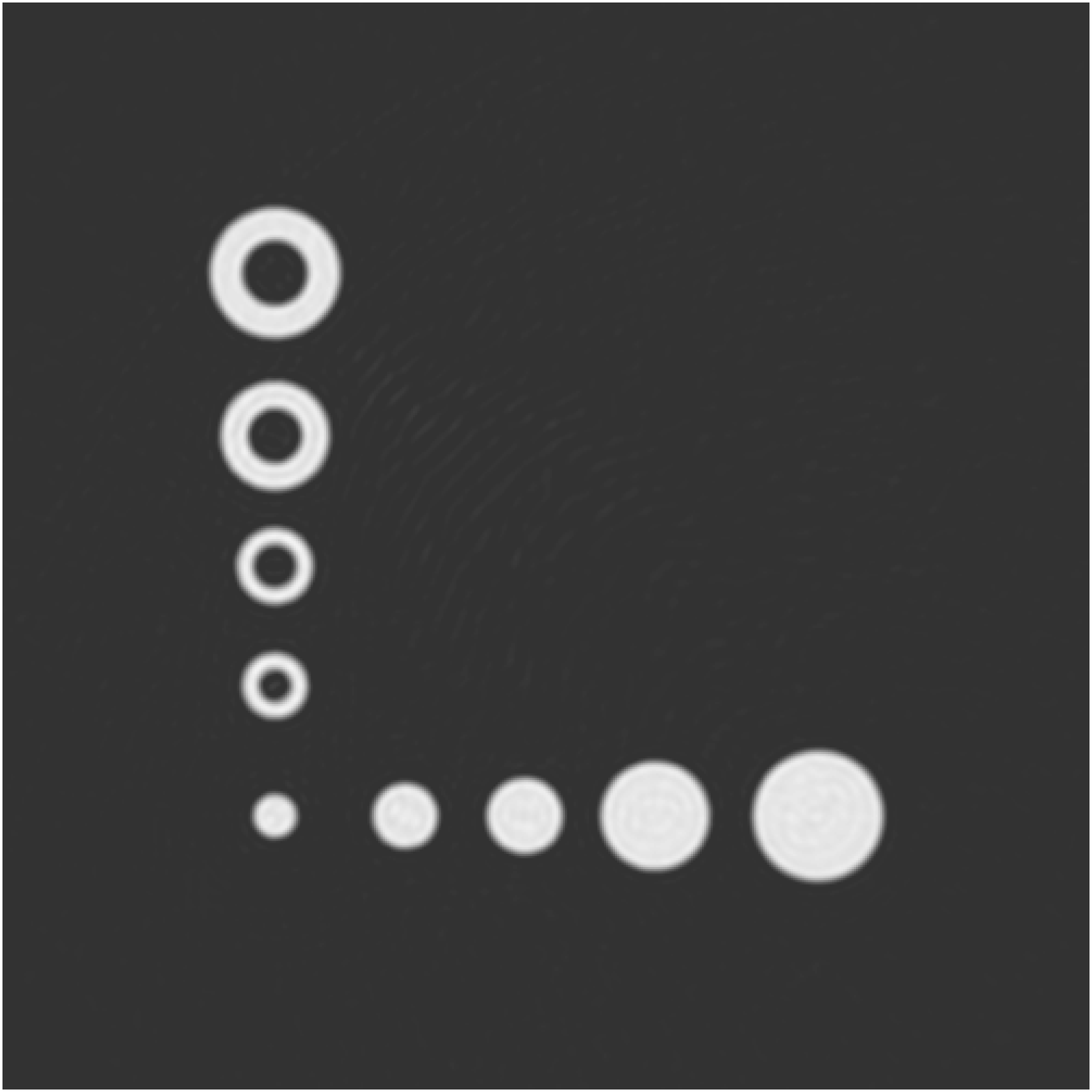}}
\subfigure{\includegraphics[width=1.8in,height=1.8in]{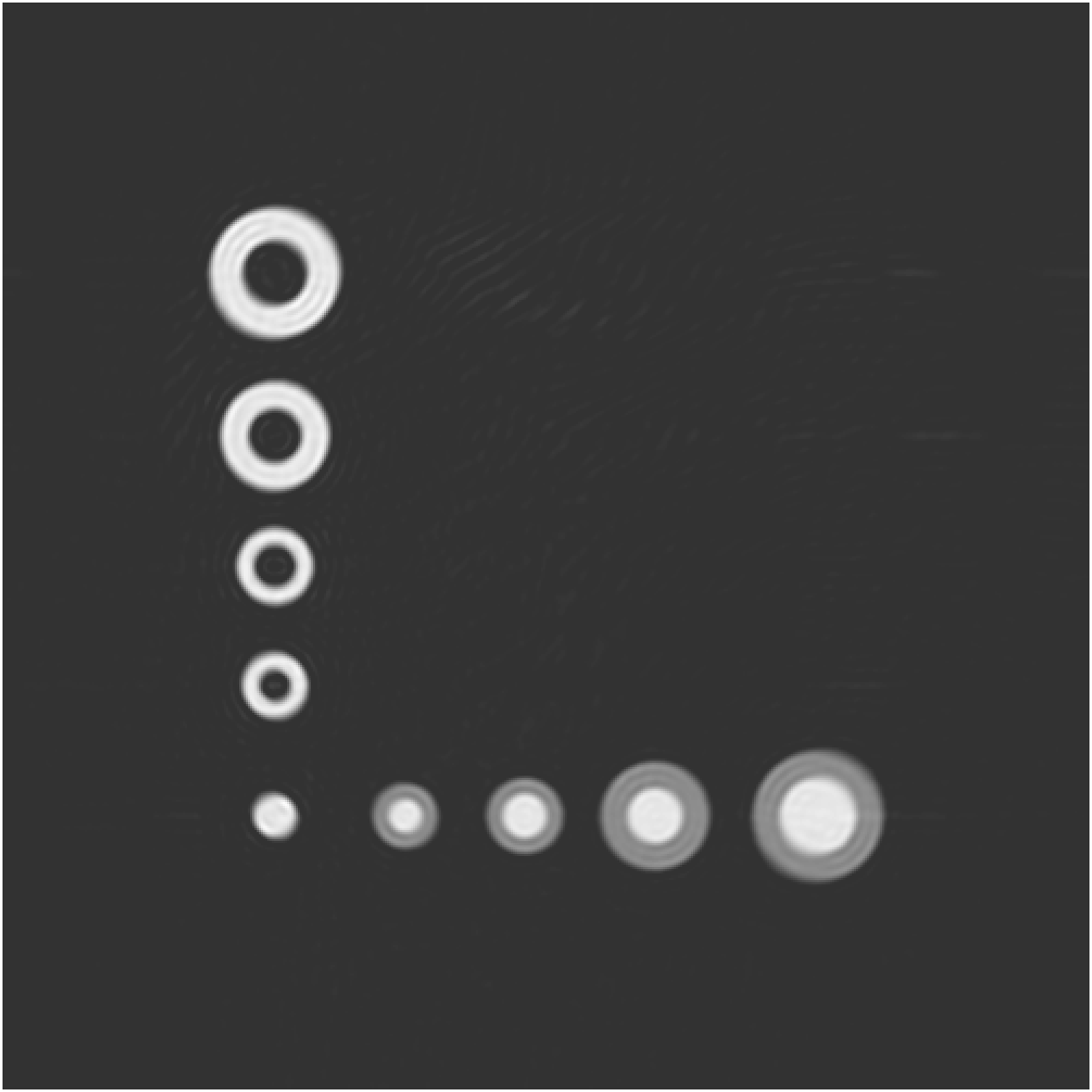}}\newline%
\setcounter{subfigure}{0}
\subfigure[Plane $x_3=-0.5$]{\includegraphics[width=1.8in,height=1.8in]{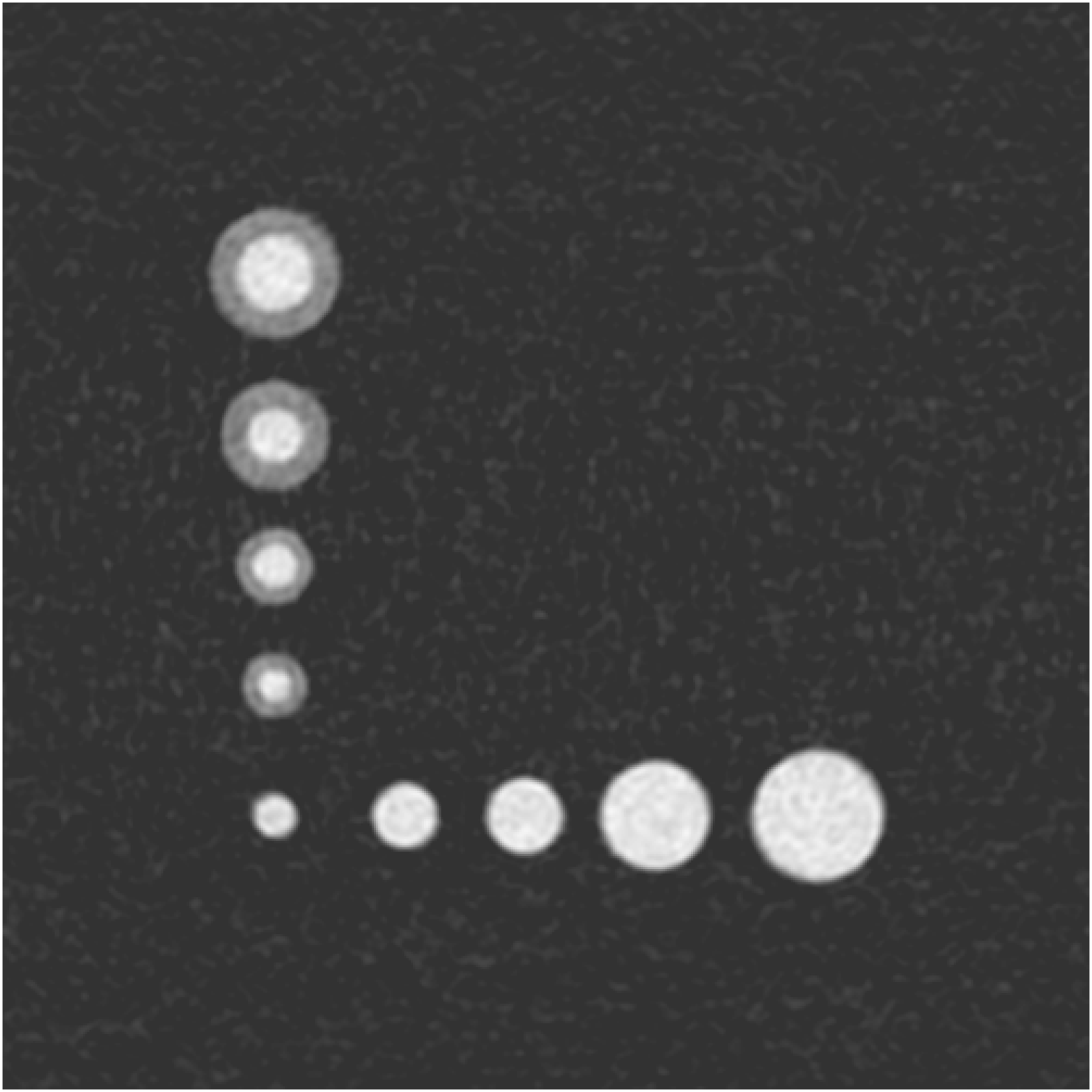}}
\subfigure[Plane $x_2=-0.5$]{\includegraphics[width=1.8in,height=1.8in]{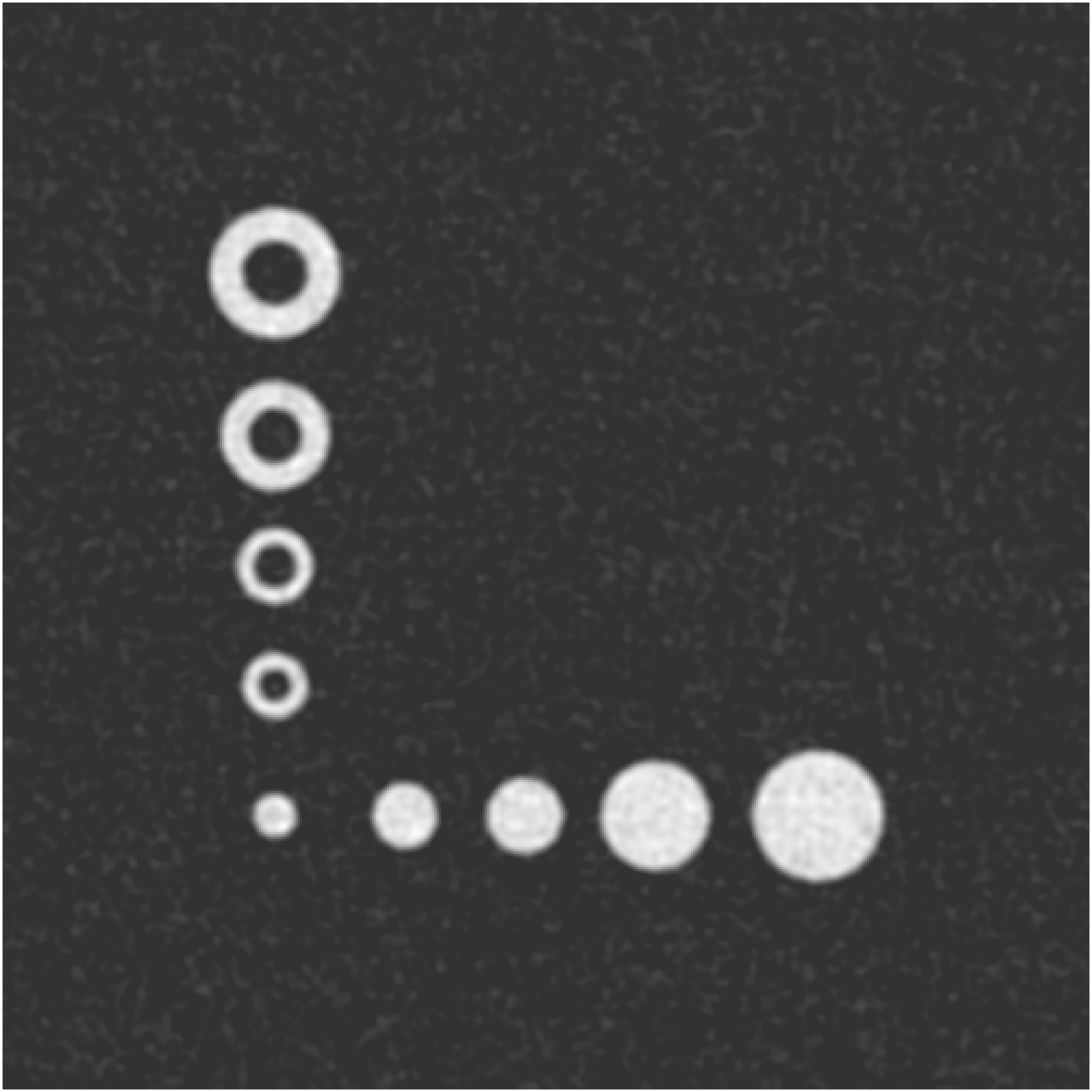}}
\subfigure[Plane $x_1=-0.5$]{\includegraphics[width=1.8in,height=1.8in]{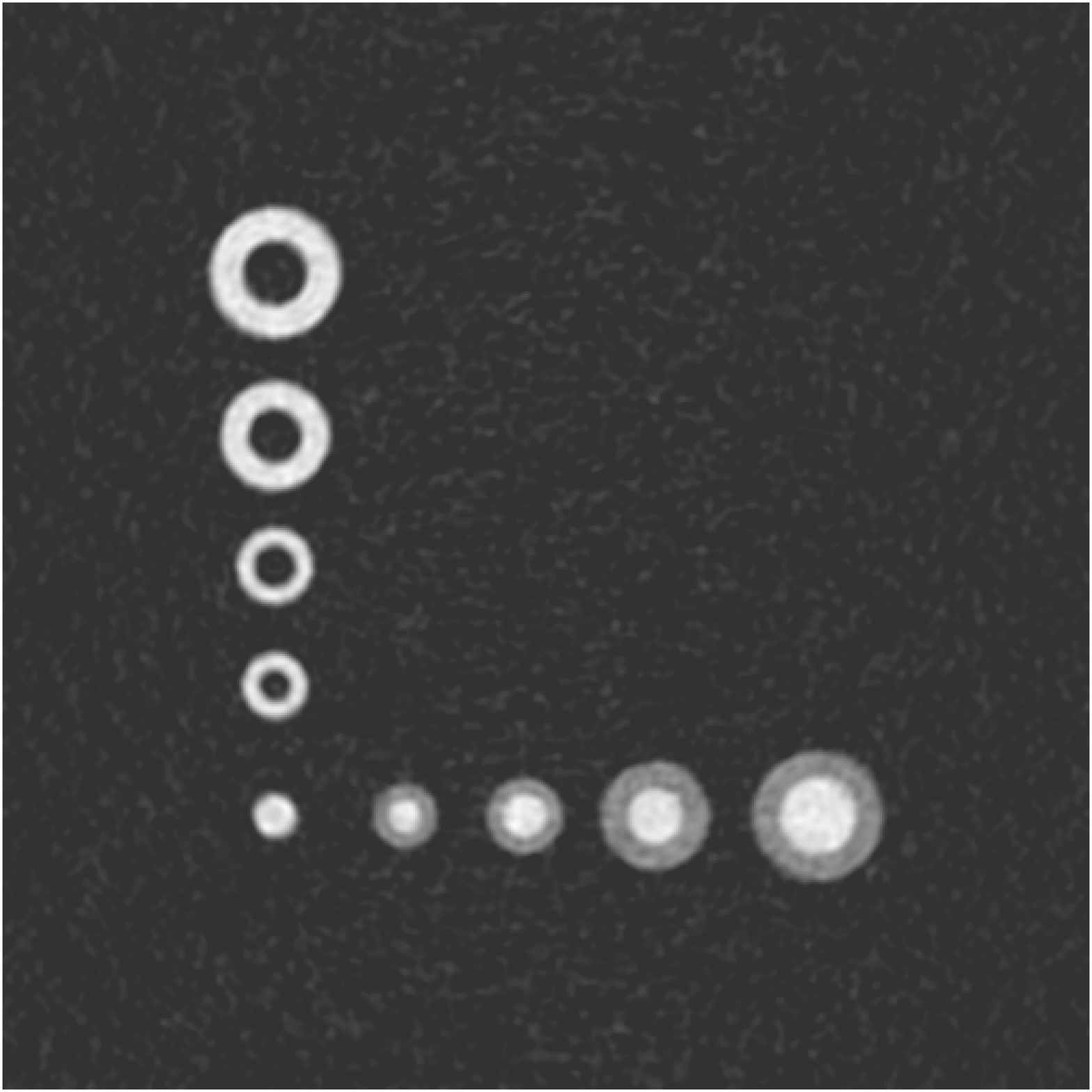}}\newline
\end{center}
\caption{Simulation in $3D$; first row represents the phantom; second row is
the reconstruction from noiseless data; third row shows reconstruction from
the data with added $50\%$ noise }%
\label{F:rec3Dflat}%
\end{figure}

It is well-known (see e.g. \cite{Vladimirov}) that integrals $v_{\alpha}(h,t)$
satisfy the wave equation in $\mathbb{R}^{2}$ (the "method of descent"\ is
based on this fact):%
\[
\left(  \frac{\partial^{2}}{\partial h_{1}^{2}}+\frac{\partial^{2}}{\partial
h_{2}^{2}}-\frac{\partial^{2}}{\partial t^{2}}\right)  v_{\alpha}=0,\quad
h\in\mathbb{R}^{2},\quad t\in\lbrack0,\infty)
\]
with the initial conditions%
\begin{align*}
v_{\alpha}(h,0)  &  \equiv M_{\alpha}(h)=\int\limits_{\mathbb{R}^{1}}%
u(h_{1}N(\alpha)+h_{2}e_{2}+sD(\alpha),0)ds\\
&  =\int\limits_{\mathbb{R}^{1}}f(h_{1}N(\alpha)+h_{2}e_{2}+sD(\alpha))ds,\\
\frac{\partial}{\partial t}v_{\alpha}(h,0)  &  =0.
\end{align*}
In other words, the initial values $M_{\alpha}$ of $v_{\alpha}$ are the
integrals of $f(x)$ along the straight lines parallel to the vector
$D(\alpha)$.
On the other hand, integrals $P_{\alpha}(y(\beta),t)$ measured by
the integrating line detectors (see equation (\ref{E:cylinder}) are equal to
the values of $v_{\alpha}(h,t)$ for all $h\ $lying on the circle $\partial B$
of radius $R$ centered at the origin $h=0$:%
\[
v_{\alpha}(y,t)=P_{\alpha}(y,t),\quad y\in\partial B,\quad t\in\lbrack
0,\infty).
\]
Therefore, for each fixed angle $\alpha,$ values of $M_{\alpha}(h)$ can be
recovered within the disk $B$ (bounded by $\partial B$) by applying the 2D
algorithm of Section~\ref{S:2Dconventional} to $P_{\alpha}(y,t)$. Further, for
a fixed value of the vertical coordinate $h_{2},$ the values of $M_{\alpha
}((h_{1},h_{2}))$ correspond to the 2D Radon transform of the restriction of
$f(x)$ to the plane $x_{2}=h_{2}.$ Thus $\left.  f(x)\right\vert _{x_{2}%
=h_{2}}$ can be reconstructed from $M_{\alpha}((h_{1},h_{2}))$ using one of
the well-known inversion algorithms for the latter transform (see e.g.
\cite{Natterer}). In particular, in ~\cite{Palt-Machzend,burg-FBP} the
well-known filtration/backprojection algorithm is utilized for such a computation.

The above mentioned techniques are not as fast as we would like, however. In
order to accelerate the computation, one can combine the fast method proposed
in Section~\ref{S:2Dconventional}, with the fast Fourier transform-based
inversion of the Radon transform \cite{Natterer}. Notice, however, that such
two stage approach adds to the algorithm a numerical inversion of the 2D Radon
transform, which may increase the computational error, since such inversion is
a (mildly) ill-posed problem. Moreover, such an algorithm is somewhat
redundant. Indeed, on the first stage the values of $M_{\alpha}(h)$ are
reconstructed by means of the Fourier synthesis (step 7 of the algorithm
presented in Section \ref{S:alg2d}). On the second stage, as a part of the
Fourier reconstruction from projections, the 1D Fourier transform in $h_{1}$
of projections $M_{\alpha}((h_{1},h_{2}))$ is computed for each value of
$\alpha$

This redundancy can be eliminated as follows. For each fixed $\alpha$ let us
apply the first 5 steps of the algorithm from Section \ref{S:alg2d} to the
data $P_{\alpha}(y,t)$ and compute functions $\hat{f}_{\alpha,2D}(\Lambda)$,
where the subscript indicates the dependence of $\hat{f}$ on the angle
$\alpha,$ and the fact that $\hat{f}_{\alpha,2D}(\Lambda)$ are solutions of 2D
problems. For a fixed $\alpha,$ function $\hat{f}_{\alpha,2D}(\Lambda)$ is the
2D Fourier transform of $M_{\alpha}(h).$ Since $M_{\alpha}(h)$ are the X-ray
projections of the initial condition $f(x),x\in\mathbb{R}^{3}$ in the
direction $D(\alpha),$ by the well-known slice-projection theorem
\cite{Natterer} $\hat{f}_{\alpha,2D}(\Lambda)$ coincides with the values of
the 3D Fourier transform $\hat{f}_{3D}(\Lambda)$ of $f(x)$ restricted to the
plane normal to $D(\alpha)$ and passing through the origin $\Lambda=0.$
Therefore, $\hat{f}_{3D}(\Lambda)$ can be reconstructed directly from the
values of $\hat{f}_{\alpha,2D}(\Lambda)$ computed for all $\alpha$ from $0$ to
$\pi.$ The sought function $f(x)$ is then obtained by computing the 3D inverse
Fourier of $\hat{f}_{3D}(\Lambda).$ This can be summarized in the form of the following

\textbf{Algorithm for solving Problem 2}

\begin{enumerate}
\item For all $\alpha$ from $0$ to $\pi$ compute $\hat{f}_{\alpha,2D}%
(\Lambda)$ from $P_{\alpha}(y,t)$ by performing the first 5 steps of the
algorithm from Section \ref{S:alg2d}

\item Interpolate values $\hat{f}_{\alpha,2D}(\Lambda)$ to obtain $\hat
{f}_{3D}(\Lambda)$ on a 3D Cartesian grid

\item Compute $f(x)$ from $\hat{f}_{3D}(\Lambda)$ by using the 3D inverse FFT.
\end{enumerate}

If the number of the cylinder directions is of order of $\mathcal{O}(n),$ step
1 of the above algorithm requires $\mathcal{O}(n^{3}\log n)$ flops. The second
step involves interpolation in the spectral domain from spherical grid to the
Cartesian grid. In our simulations we utilized cubic interpolation in the
radial direction and bilinear interpolation in the angular directions; the
number of operations involved in this step is $\mathcal{O}(n^{3}).$ The third,
final step requires $\mathcal{O}(n^{3}\log n)$ operations; the whole algorithm
then needs only $\mathcal{O}(n^{3}\log n)$ flops.

In order to evaluate the performance of the present algorithm we simulated the
measurements corresponding to 272 integrating line detectors rotated over 512
directions $\alpha$ equispaced from $0$ to $\pi.$ As a phantom we used a
collection of characteristic functions of balls centered on the (pair-wise)
intersections of the planes $x_{1}=-0.5,$ $x_{2}=-0.5$, $x_{3}=-0.5$, and
lying within the unit sphere, as shown in Figure~\ref{F:3Dphan}. (The
orientation of the axes corresponds to vectors $e_{1}$, $e_{2}$ and $e_{3}$ as
illustrated in Figure~1(b).) The cross-sections of the phantom are shown in
the first row of Figure~\ref{F:rec3Dflat}.

For each direction $\alpha$, 500 time samples were simulated in the interval
$t\in\lbrack0,5]$; the rest of the signal was neglected. The image was
reconstructed on a $500\times500\times500$ Cartesian grid containing
125 million unknowns (although only about a half of them lied within the
unit sphere where function $f(x)$ was supported).
The cross sections of the reconstructed $f(x)$ are shown in the second row
of Figure~\ref{F:rec3Dflat}.

The above computation took 67 seconds on the desktop computer described
in Section \ref{S:alg2d}. A comparison can be made with a reconstruction
obtained by the time reversal using  finite differences  in a cubic domain
with $251 \times 251 \times 251$ computational grid.
On our computer it took about $50$ min. Since
such a time reversal method scales as $\mathcal{O}(n^4)$,
on a grid of the size $500 \times 500 \times 500$ the reconstruction would take about
$13$ hours, or three orders of magnitude longer than the time required
by the present algorithm. (Such a comparison is quite crude since
a different problem is solved in a different computational domain;
nevertheless it indicates that our method is indeed very fast.)

The third row of Figure~\ref{F:rec3Dflat} demonstrates images reconstructed from the data
with added simulated noise with intensity 50\% (in $L^{2}$ norm) of the
signal. The level of noise in the reconstructed images is surprisingly low.
The explanation of such stability is the same as in Section \ref{S:alg2d}.
\begin{figure}[t]
\begin{center}
\subfigure[]{\includegraphics[width=1.2in,height=2.5in]{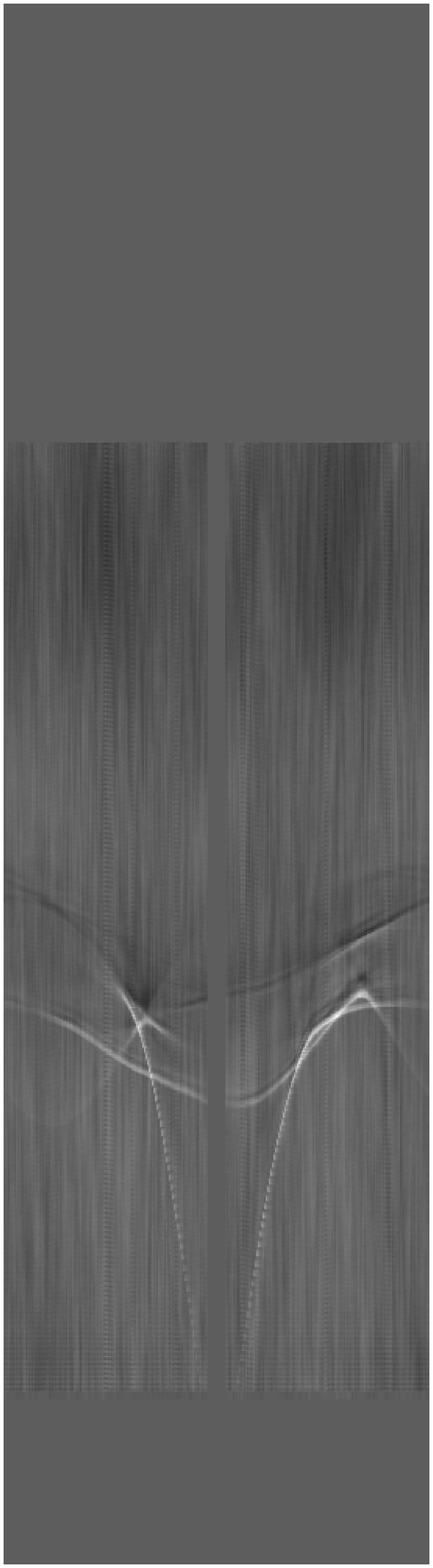}}
\subfigure[]{\includegraphics[width=2.9in,height=2.5in]{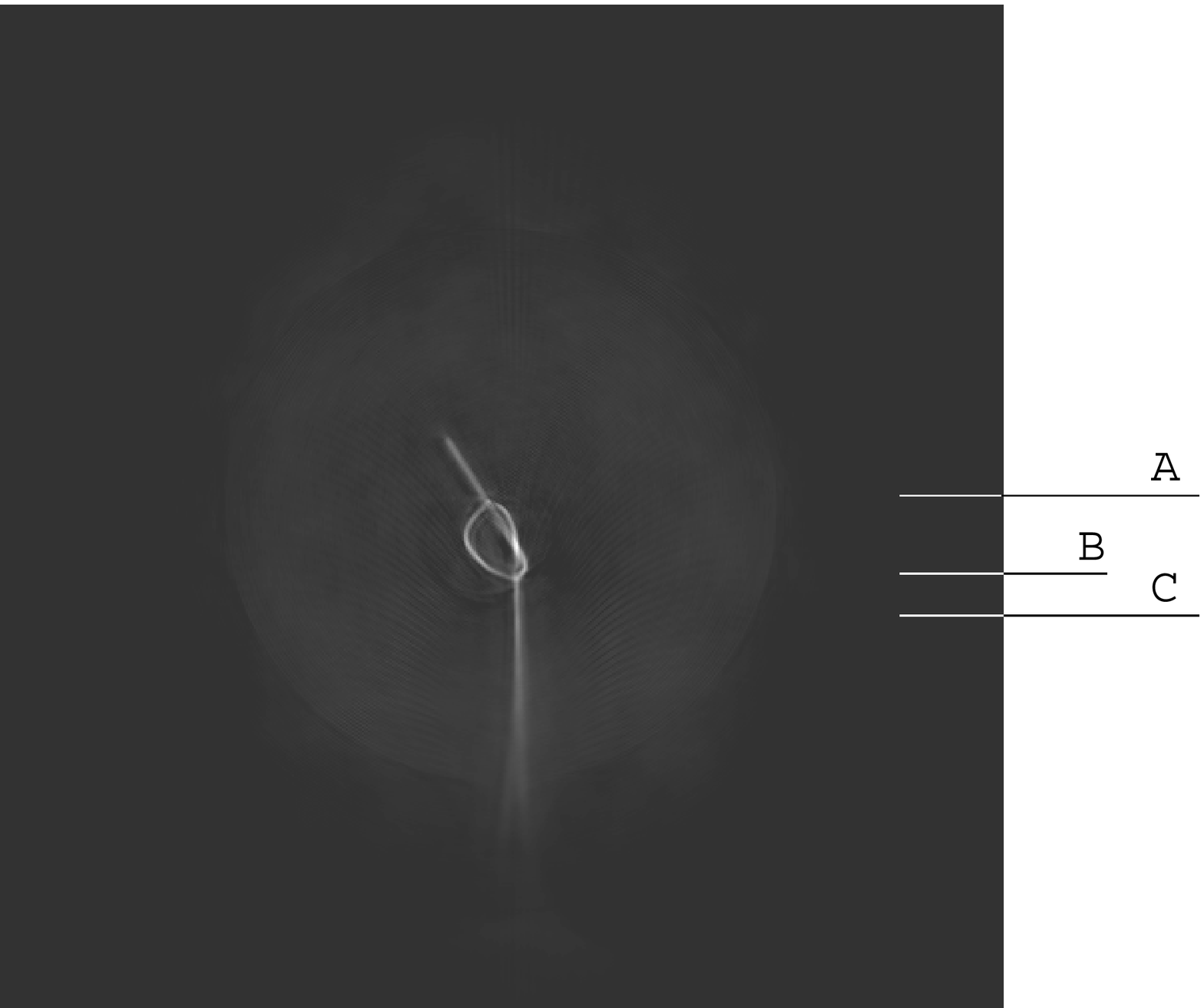}}
\end{center}
\caption{ Reconstruction of a 2D projection from the data measured by
integrating line detectors (a) the data (b) the reconstruction}%
\label{F:real2D}%
\end{figure}

\section{Applications to real data}
\label{S:real}

In this section we illustrate the work of our algorithms by processing a set
of real data kindly provided to us by RECENDT. This set of data was measured
by integrating line detectors, as shown in Figure~1 and described in detail in
Section~\ref{S:linedet}. As the test object the researchers from RECENDT used
a piece of a human hair tied in a knot. The number of the detector directions
$D(\alpha)$ in this set was $25.$ For each direction, the linear detector was
respectively placed in each of 272 equispaced positions on the surface of the
(imaginary) cylinder. (In fact, 11 out of 272 positions at the bottom were unavailable and
the corresponding data were replaced by zeros.) For each position 10000 time
samples were measured. To reduce the noise, the signal was smoothed by a
convolution with a Gaussian and downsampled by a factor of 10, so that the
number of time samples on the input of the reconstruction algorithm was 1000.
Moreover, to reduce strong artifacts at the beginning and at the end of each
time series, the signal was set to zero there as well. The resulting set of
data (corresponding to the detectors being aligned along $e_{3})$ is
represented by the gray-scale image shown in Figure~\ref{F:real2D}(a). Each
vertical line in this image corresponds to the time series for one detector
position of a line detector; the bottom corresponds to $t=0.$
\begin{figure}[t]
\begin{center}
\subfigure{\includegraphics[width=1.8in,height=1.8in]{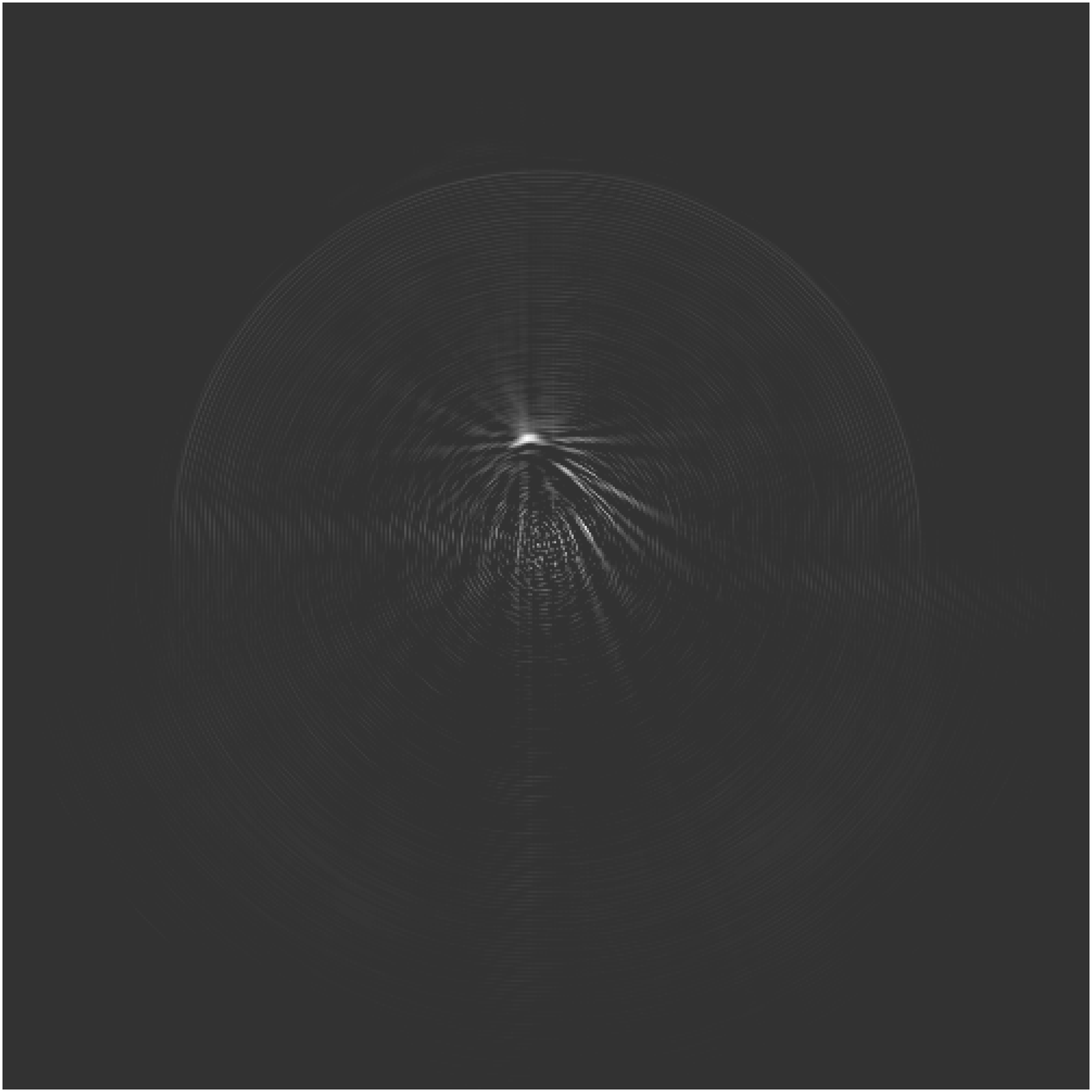}}
\subfigure{\includegraphics[width=1.8in,height=1.8in]{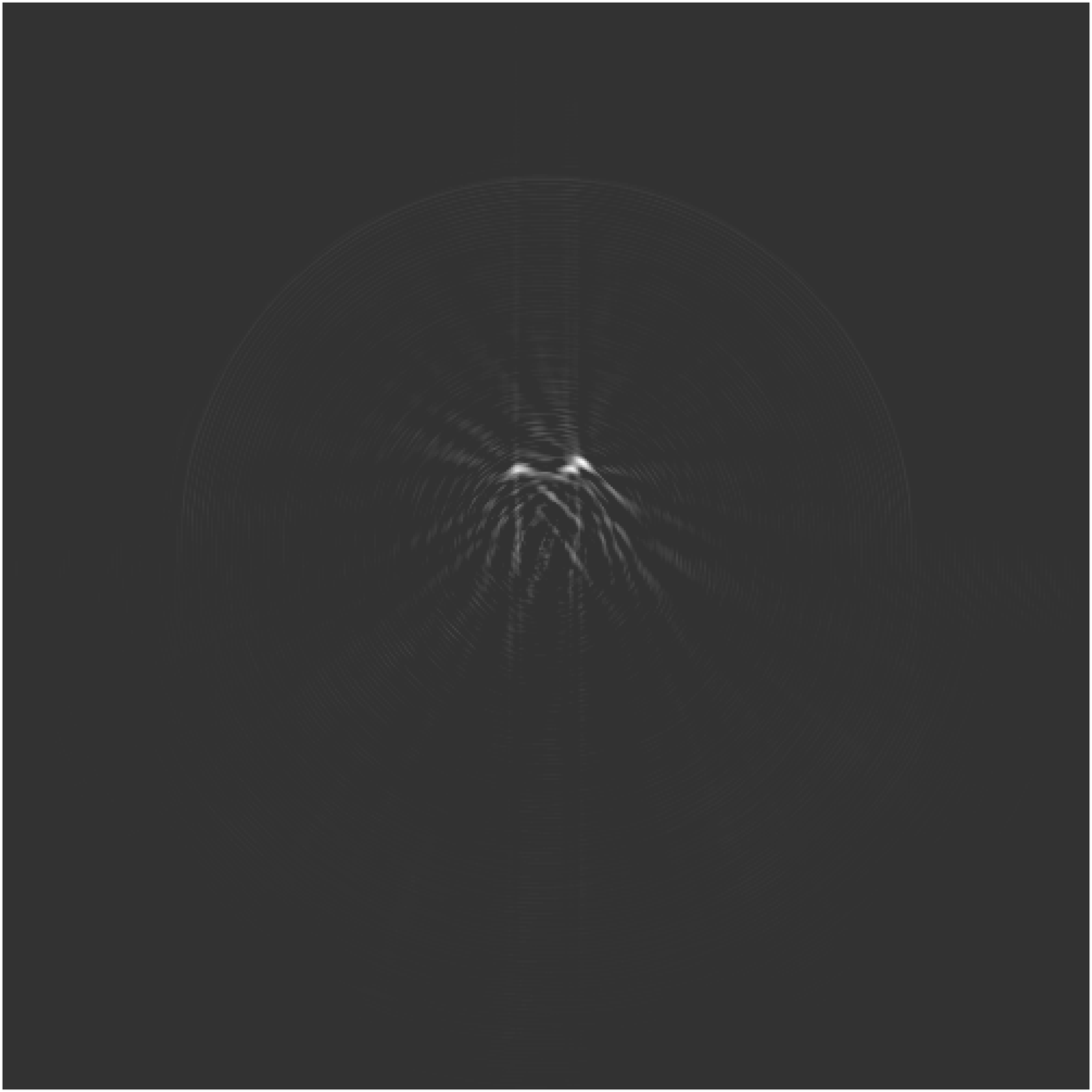}}
\subfigure{\includegraphics[width=1.8in,height=1.8in]{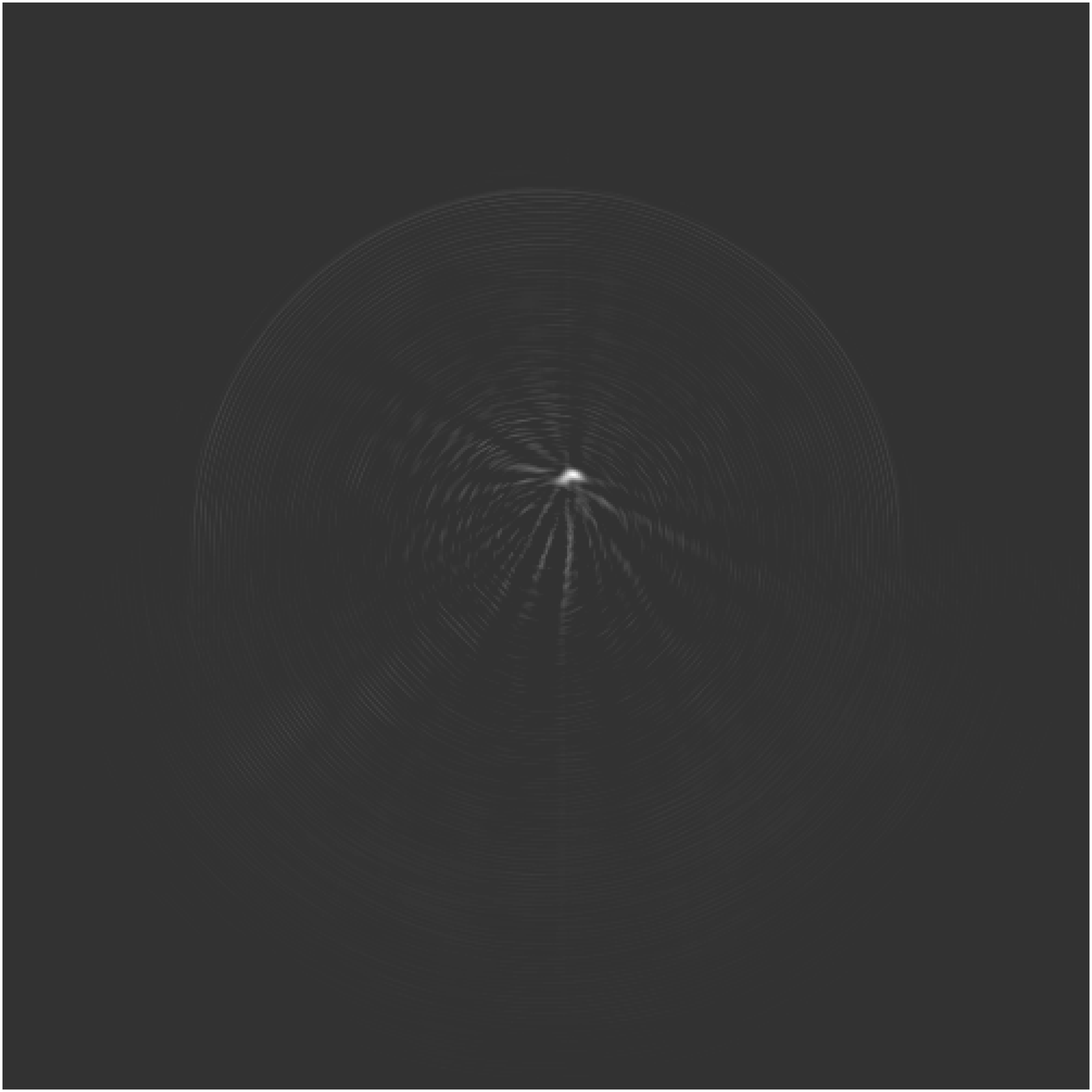}}
\end{center}
\caption{3D reconstruction from real data; slices parallel to the plane $O
x_{1} x_{3}$ at the $x_{2}$ levels that correspond, respectively, to marks
A,B, and C in Figure~\ref{F:real2D} }%
\label{F:real3D}%
\end{figure}

Figure~\ref{F:real2D}(b) demonstrates the image reconstructed by applying
algorithm of Section~\ref{S:alg2d} to the data shown in Figure~\ref{F:real2D}%
(a). As explained in Section~\ref{S:linedet}, the result is not the complete
reconstruction, but the X-ray transform of the 3D density function $f(x)$
corresponding to the hair.

In order to reconstruct $f(x)$ we applied algorithm of Section~\ref{S:linedet}
to the full set of data. The results are shown in Figure~\ref{F:real3D}. The
three images in this figure correspond to the horizontal cross-sections of the
test object at the levels A, B, and C, respectively. One can see sharp spikes
in the images corresponding to the location(s) of the hair. There are
significant radial artifacts in these images, arising due to the insufficient
sampling in angle $\alpha$ (i.e. insufficient number of the detector
directions). It is known \cite{Natterer} that for the inversion of the 2D
Radon transform the optimal number of the angular directions is of the same
order of magnitude as the desired resolution of the spatial grid (several
hundred, in out case). Unfortunately, our data set only contained data
corresponding to 25 detector directions. As the reconstructions obtained in
Section~\ref{S:linedet} from simulated data show,  if the number of directions
equals several hundred, the algorithm yields very accurate and detailed images.

\section*{Acknowledgements}

The author would like to thank RECENDT for providing the data of real
measurements. We are grateful to Drs. P.~Burgholzer, H.~Gr\"{u}n, and
H.~Roitner for the fruitful discussions of the data acquisition schemes with
integrating line detectors, related data processing techniques and open
problems in this area. The author gratefully acknowledges support by the NSF
through the grant DMS-090824.


\begin{thebibliography}{99}                                                                                               %


\bibitem {AK}M. Agranovsky and P. Kuchment, \emph{Uniqueness of reconstruction and
an inversion procedure for thermoacoustic and photoacoustic tomography with
variable sound speed}, Inverse Problems, \textbf{23}:2089--102, 2007.

\bibitem {AmbKuch}G. Ambartsoumian and P. Kuchment, \emph{A range description for
the planar circular Radon transform} SIAM J. Math. Anal., \textbf{38}%
(2): 681--92, 2006.

\bibitem {Andersson}L. E. Andersson, \emph{On the determination of a function from
spherical averages} SIAM J. Math. Anal., \textbf{19}(1): 214--32, 1988.

\bibitem {burg-FBP}P. Burgholzer, J. Bauer-Marschallinger, H. Gr\"{u}n, M.
Haltmeier, and G. Paltauf, \emph{Temporal back-projection algorithms for
photoacoustic tomography with integrating line detectors},  Inverse
Problems, \textbf{23}:S65-S80, 2007.

\bibitem {burg-area-line}P. Burgholzer, C. Hofer, G. Paltauf, M. Haltmeier, O.
Scherzer, \emph{Thermoacoustic tomography with integrating area and line detectors},
IEEE Transactions on Ultrasonics, Ferroelectrics, and Frequency Control,
\textbf{52}(9):1577--83, 2005.

\bibitem {burg-fabry}P. Burgholzer, C. Hofer, G. J. Matt, G. Paltauf, M.
Haltmeier, and O. Scherzer, \emph{Thermoacoustic tomography using a fiber-based
Fabry-Perot interferometer as an integrating line detector}, Proc. SPIE {6086},
434--42, 2006.

\bibitem {burg-exac-appro}P. Burgholzer, G. J. Matt, M. Haltmeier, and G.
Paltauf, \emph{Exact and approximative imaging methods for photoacoustic tomography
using an arbitrary detection surface}, Phys Review E, \textbf{75},
046706, 2007.

\bibitem {Colton}D. Colton and R. Kress, Inverse acoustic and
electromagnetic scattering theory, Springer-Verlag, 1992

\bibitem {Dutt}A. Dutt and V. Rokhlin, \emph{Fast Fourier Transforms For
Nonequispaced Data}, Siam J. Sci. Comput., \textbf{14}(6): 1368--93, 1993.

\bibitem {Fawcett}J. A. Fawcett, \emph{Inversion of $n$-dimensional spherical
averages}, SIAM J. Appl. Math., \textbf{45}(2): 336--41, 1985.

\bibitem {Finch-even}D. Finch, M. Haltmeier and Rakesh,
\emph{Inversion of spherical
means and the wave equation in even dimensions},
SIAM J. Appl. Math., \textbf{ 68}(2): 392--412, 2007.

\bibitem {FPR}D. Finch, S. Patch and Rakesh,
\emph{Determining a function from its
mean values over a family of spheres},
SIAM J. Math. Anal.,
\textbf{35}(5): 1213--40, 2004.

\bibitem {Haltm-circ}M. Haltmeier, O. Scherzer, P. Burgholzer, R. Nuster, and
G. Paltauf,
\emph{Thermoacoustic Tomography And The Circular Radon Transform: Exact
Inversion Formula}, Mathematical Models and Methods in Applied Sciences,
\textbf{17}(4): 635--55, 2007.

\bibitem {Haltm-fft}M. Haltmeier, O. Scherzer and G. Zangerl,
\emph{A Reconstruction
Algorithm for Photoacoustic Imaging Based on the Nonuniform FFT},
IEEE Trans. Med. Imag., \textbf{28}(11):1727--35, 2009.

\bibitem {Healy}D. M. Healy Jr., D. N. Rockmore, P. J. Kostelec, and S. Moore.
\emph{FFTs for the 2-Sphere -- Improvements and Variations},
J. Fourier Anal. and Appl., \textbf{9}(4): 341--85,\ 2003.

\bibitem {HKN}Y. Hristova, P. Kuchment, and L. Nguyen,
\emph{On reconstruction and
time reversal in thermoacoustic tomography in homogeneous and non-homogeneous
acoustic media}, Inverse Problems, \textbf{24}: 055006, 2008.

%\bibitem {Knock}L. Knockaert. Fast Hankel transform by fast sine and cosine
%transforms: the Mellin connection. , IEEE Trans. Signal Process.}
%\textbf{48}(6): 1695--701, 2000.


\bibitem {Kruger}R. A. Kruger, P. Liu, Y. R. Fang, and C. R. Appledorn,
\emph{Photoacoustic ultrasound (PAUS) reconstruction tomography},
Med. Phys., \textbf{22}: 1605--09, 1995.

\bibitem {Springer}P. Kuchment and L. Kunyansky,
Mathematics of Photoacoustic and Thermoacoustic Tomography,
in Handbook of Mathematical Methods in Imaging, Scherzer, Otmar (Ed.) Springer, 2011.

\bibitem {Kunyansky}L. Kunyansky,
\emph{Explicit inversion formulae for the
spherical mean Radon transform},
Inverse Problems, \textbf{23}: 737-783, 2007.

\bibitem {Kun-series}L. Kunyansky,
\emph{A series solution and a fast algorithm for
the inversion of the spherical mean Radon transform},
Inverse Problems, \textbf{23}: S11--S20, 2007.

\bibitem {Mohlen}M. J. Mohlenkamp, \emph{A Fast Transform for Spherical Harmonics},
J. Fourier Anal. Appl., \textbf{2}: 159--84, 1999.

\bibitem {Natterer} F. Natterer,
The mathematics of computerized tomography, New York, Wiley,  1986.

\bibitem {Nat-Wub}F. Natterer and F. W\"{u}bbeling, Mathematical Methods
in Image Reconstruction, ser. Monographs Math. Model. Comput., Philadelphia,
PA: SIAM, 2001, vol. 5.

\bibitem {nguyen}L. Nguyen, \emph{A family of inversion formulas in thermoacoustic
tomography},
Inverse Problems and Imaging, \textbf{3}(4): 649-675, 2009.

\bibitem {Norton1}S. J. Norton, \emph{Reconstruction of a two-dimensional reflecting
medium over a circular domain: exact solution},
J. Acoust. Soc. Am.,
\textbf{67}: 1266-1273, 1980.

\bibitem {Norton2}S. J. Norton and M. Linzer,
 \emph{Ultrasonic reflectivity imaging
in three dimensions: exact inverse scattering solutions for plane,
cylindrical, and spherical apertures}, IEEE Transactions on Biomedical
Engineering, \textbf{28}: 200-202, 1981.

\bibitem {Oraev94}A. A. Oraevsky, S. L. Jacques, R. O. Esenaliev, and F. K.
Tittel,
\emph{Laser-based ptoacoustic imaging in biological tissues},
Proc. SPIE, \textbf{2134A}:122-128, 1994.

\bibitem {Palt-Machzend}G. Paltauf, R. Nuster, M. Haltmeier, and P.
Burgholzer,
\emph{Thermoacoustic Computed Tomography using a Mach-Zehnder
interferometer as acoustic line detector},
Appl. Opt., \textbf{46}(16):3352-8, 2007.

\bibitem {Paltaufexp}G. Paltauf, R. Nuster, M. Haltmeier and P. Burgholzer,
\emph{Experimental evaluation of reconstruction algorithms for limited view
photoacoustic tomography with line detectors}, Inverse Problems,
\textbf{23}: S81--S94, 2007.

\bibitem {Paltaufwei}G. Paltauf, R. Nuster, and P. Burgholzer,
\emph{Weight factors for limited angle photoacoustic tomography},
Phys. Med. Biol.,
\textbf{54}: 3303--14, 2009.

\bibitem {Potts}D. Potts, G. Steidl and M. Tasche,
\emph{Fast and stable algorithms
for discrete spherical Fourier transforms.}, Linear Algebra Appl.,
\textbf{275/276}: 433--450, 1998.

\bibitem {Ramm}A. G. Ramm, \emph{Injectivity of the spherical means operator},
C. R. Math. Acad. Sci. Paris, \textbf{335(}12): 1033--38, 2002.

\bibitem {Suda}R. Suda and M. Takami,
\emph{A fast spherical harmonics transform
algorithm}, Mathematics of computation, \textbf{71}(238): 703--15, 2001.

%\bibitem {Talman}J. P. Talman. Numerical Fourier and Bessel transforms in
%logarithmic variables. , J. Comp. Phys.} \textbf{29}: 35--48, 1978.


\bibitem {Vladimirov}V. S. Vladimirov, Equations of mathematical
physics. (Translated from the Russian by Audrey Littlewood. Edited by Alan
Jeffrey.) Pure and Applied Mathematics, \textbf{3} Marcel Dekker, New York, 1971.

\bibitem {CRC}L. Wang, (Editor), Photoacoustic imaging and spectroscopy,
CRC Press, Boca Raton, FL, 2009.

\bibitem {Wang-book}L. V. Wang and H. Wu, Biomedical Optics. Principles
and Imaging, Wiley-Interscience, 2007.

\bibitem {MXW1}M. Xu and L.-H. V. Wang,
\emph{Time-domain reconstruction for
thermoacoustic tomography in a spherical geometry}, IEEE Trans. Med.
Imag., \textbf{21}: 814-822, 2002.

\bibitem {MXW2}M. Xu and L.-H. V. Wang, \emph{Universal back-projection algorithm
for photoacoustic computed tomography}, Phys. Rev. E, \textbf{71}:016706, 2005.
\end{thebibliography}
\end{document}